\documentclass[journal]{IEEEtran}

\makeatletter
\let\origsubsubsection\subsubsection
\renewcommand\subsubsection{\@ifstar{\starsubsubsection}{\nostarsubsubsection}}

\newcommand\nostarsubsubsection[1]
{\origsubsubsection{#1}}

\newcommand\starsubsubsection[1]
{\subsubsectionprelude\origsubsubsection*{#1}\subsubsectionpostlude}

\newcommand\subsubsectionprelude{%
  \vspace{0em}
}

\newcommand\subsubsectionpostlude{%
  \vspace{0.5em}
}
\makeatother

\usepackage{enumitem}

\usepackage{cite}

\usepackage{graphicx}

\usepackage{color}

\usepackage{amsmath,amsfonts,amssymb,amsthm}
\usepackage{mathtools}
\usepackage{units}

\newcommand{\overbar}[1]{\mkern 1.5mu\overline{\mkern-1.5mu#1\mkern-1.5mu}\mkern 1.5mu}

\newcommand*{\VEC}[1]  {\ensuremath{\boldsymbol{#1}}}
\newcommand*{\MAT}[1]  {\ensuremath{\boldsymbol{#1}}}

\newcommand{\argmin}{\operatornamewithlimits{argmin}}

\newcommand{\reel}{\mathfrak{Re}}

\DeclareMathOperator{\tr}{tr}
\DeclareMathOperator{\herm}{herm}
\DeclareMathOperator{\skewh}{skew}
\DeclareMathOperator{\vect}{vec}
\DeclareMathOperator{\D}{D}
\DeclareMathOperator{\grad}{grad}
\DeclareMathOperator{\Hess}{Hess}
\DeclareMathOperator{\uf}{uf}
\DeclareMathOperator{\err}{err}
\DeclareMathOperator{\spann}{span}

\makeatletter
\def\@endtheorem{\endtrivlist}
\makeatother

\newtheorem{definition}{Definition}
\newtheorem{proposition}{Proposition}
\newtheorem{corollary}{Corollary}

\DeclarePairedDelimiter\norm{\lVert}{\rVert}
\makeatletter
\let\oldnorm\norm
\def\norm{\@ifstar{\oldnorm}{\oldnorm*}}
\makeatother

\newcommand{\I}{\mathfrak{i}}

\usepackage[font=footnotesize,skip=2pt]{caption}

\usepackage{tikz}
\usepackage{pgfplots}
\pgfplotsset{compat=1.9}
 \pgfplotsset{every  tick label/.append style={font=\scriptsize},
 			}

\usepackage{readarray}
\usepackage{filecontents}

\usepackage{changes}
\definechangesauthor[color=green]{FP}

\newlength\height 
\newlength\width

\begin{document}

\title{A Riemannian Framework for Low-Rank Structured Elliptical Models}

\author{\IEEEauthorblockN{
Florent Bouchard, 
Arnaud Breloy,~\IEEEmembership{Member,~IEEE},
Guillaume Ginolhac,~\IEEEmembership{Senior Member,~IEEE},
Alexandre Renaux,~\IEEEmembership{Member,~IEEE},
Fr\'ed\'eric Pascal,~\IEEEmembership{Senior Member,~IEEE}}
\thanks{
Florent Bouchard and Guillaume Ginolhac are with LISTIC (EA3703), University Savoie Mont Blanc, France (e-mails: florent.bouchard@univ-smb.fr, guillaume.ginolhac@univ-smb.fr).
Arnaud Breloy is with LEME (EA4416), University Paris Nanterrre, France (e-mail: abreloy@parisnanterre.fr).
Alexandre Renaux is with Laboratoire des signaux et syst\`emes (UMR8506), University Paris-Sud, France (e-mail: alexandre.renaux@u-psud.fr).
Frederic Pascal is with Universit\'e Paris-Saclay, CNRS, CentraleSup\'elec,  Laboratoire des signaux et syst\`emes (UMR8506), 91190, Gif-sur-Yvette, France (e-mail: frederic.pascal@l2s.centralesupelec.fr). This work was supported by ANR-ASTRID MARGARITA (ANR-17-ASTR-0015).
} 
}

\maketitle

\begin{abstract}
This paper proposes an original Riemmanian geometry for low-rank structured elliptical models, \emph{i.e.}, when samples are elliptically distributed with a covariance matrix that has a low-rank plus identity structure.  The considered geometry is the one induced by the product of the Stiefel manifold and the manifold of Hermitian positive definite matrices, quotiented by the unitary group. One of the main contribution is to consider an original Riemannian metric, leading to new representations of tangent spaces and geodesics. From this geometry, we derive a new Riemannian optimization framework for robust covariance estimation, which is leveraged to minimize the popular Tyler's cost function on the considered quotient manifold. We also obtain a new divergence function, which is exploited to define a geometrical error measure on the quotient, and the corresponding intrinsic Cram\'er-Rao lower bound is derived. Thanks to the structure of the chosen parametrization, we further consider the subspace estimation error on the Grassmann manifold and provide its intrinsic Cram\'er-Rao lower bound. Our theoretical results are illustrated on some numerical experiments, showing the interest of the proposed optimization framework and that performance bounds can be reached.
%
\end{abstract}


%
\IEEEpeerreviewmaketitle

\section{Introduction}
\label{sec:intro}
\IEEEPARstart{C}{omplex elliptically symmetric} distributions offer a general family of statistical models that encompasses most of standard multivariate distributions, including the Gaussian one, as well as many heavy-tailed distributions, such as multivariate Student $t$-, and $K$- distributions (cf.~\cite{OTKP12a} for a review on this topic).
These models have been leveraged successfully in numerous applications thanks to their good empirical fit to datasets, \emph{e.g.}, in image processing~\cite{PSWS03,SS07,ZV10} or array processing~\cite{GGR06,OTKP12b}.
On top of that, elliptical models have also attracted a lot of interest, as they allow robust estimation processes to be derived.
For example, $M$-estimators~\cite{MY76, T87}, defined as generalized maximum likelihood estimators of elliptical models, have been shown to be robust to model mismatches and contaminated data (outliers)~\cite{OTKP12a}.
While alleviating robustness issues, the development of estimation algorithms under elliptical models is still challenged by ``small $n$ large $p$'' problems (where $n$ and $p$ respectively stand for the sample size and the dimension).

In several applications, one can rightfully assume that the relevant information lies in a low dimensional subspace.
This is reflected by a low-rank structure of the covariance matrix, often referred to as spiked model~\cite{J01}. 
This idea plays a central role in principal component analysis~\cite{TB99}, subspace recovery~\cite{LM18}, and related dimension reduction algorithms. 
Low-rank models also play a central role in array processing~\cite{HPREK14} and financial time series analysis~\cite{R11} (where they are also referred to as factor models).

Estimation processes in low-rank models have been well studied for Gaussian distributions~\cite{TB99,KMR14}.
Unfortunately, the results obtained in this case cannot be trivially transposed to elliptical distributions.
For example, low-rank structured counterparts of $M$-estimators are not expressed in closed form, nor directly tractable.
Additionally, ultimate statistical performance characterization is not obvious in this context, due to constraints/ambiguities on the parameters space.

This paper proposes to leverage tools from Riemannian geometry in order to answer the previous questions with a unified view.
The Riemmanian standpoint was adopted in \cite{S05} to derive intrinsic (\emph{i.e.}, manifold oriented) Cram\'er-Rao lower bounds, then applied to study both unstructured and low-rank Gaussian models. This lead to interesting results and insights, such as performance bounds for various Riemmanian distances, and the characterization of a bias of the sample covariance matrix at low sample support, not exhibited by the traditional Euclidean analysis.
The Riemmanian geometry of the manifold of Hermitian positive definite matrices has also been recently used to study unstructured elliptical models.
It notably revealed hidden (geodesic) convexity properties of elliptical distribution's likelihood functions~\cite{W12a}, and allowed to derive new regularization-based estimation algorithms ~\cite{W12b, OT14, DT16}.
Studying low-rank elliptical models requires to turn to the manifold of Hermitian positive semi-definite matrices of fixed rank $k$ ($k<p$), which has, to the best of our knowledge, not been proposed in this context.
The contributions associated to the proposed framework for low-rank elliptical models follow three main axes, summed up below.

\subsection{Geometry for low-rank structured elliptical models}

The statistical parameter of the considered low-rank model for complex elliptically symmetric distributions lives in the manifold $\mathcal{H}^{+}_{p,k}$ of $p\times p$ Hermitian positive semi-definite matrices of rank $k$.
This manifold has recently attracted much attention and several geometries have been proposed for it; see \emph{e.g.},~\cite{BS09,VV10,MBS11,VAV12,MA18}.
In this work, we consider the geometry induced by the quotient $(\textup{St}_{p,k}\times\mathcal{H}^{++}_k)/\mathcal{U}_k$, \emph{i.e.,} the product manifold of the complex Stiefel manifold $\textup{St}_{p,k}$ of $p\times k$ orthogonal matrices (with $p>k$) and the manifold $\mathcal{H}^{++}_k$ of $k\times k$ Hermitian positive definite matrices, quotiented by the unitary group $\mathcal{U}_k$.
This geometry has already been studied in the context of low-rank matrices in~\cite{BS09,MBS11}.
It is of particular interest in our context because the principal subspace of the covariance matrix is directly obtained from this parametrization and a divergence function, which can be exploited to measure estimation errors, is available in closed form~\cite{BS09}.

Our framework differs from the works~\cite{BS09,MBS11} as we propose a new Riemannian metric on the product $\textup{St}_{p,k}\times\mathcal{H}^{++}_k$: the part on $\textup{St}_{p,k}$ is the so-called canonical metric on Stiefel~\cite{EAS98} while the part on $\mathcal{H}^{++}_k$ is a general form of the affine invariant metric which corresponds to the Fisher information metric of elliptical distributions on $\mathcal{H}^{++}_k$~\cite{BGRB18}.
As a direct consequence, the representations of tangent spaces of the quotient $(\textup{St}_{p,k}\times\mathcal{H}^{++}_k)/\mathcal{U}_k$, geodesics, Riemannian gradient and Hessian used for optimization are original in this context.
We also introduce a retraction, which corresponds to a second order approximation of the geodesics.
Moreover, we derive a new divergence function on the quotient, which is inspired by the one of~\cite{BS09}.

\subsection{Algorithms for robust low-rank covariance matrix estimation}

Covariance matrix estimation is a crucial step in many machine learning and signal processing algorithms.
In elliptical models, $M$-estimators~\cite{MY76, T87} offer a robust alternative to the traditional sample covariance matrix. 
These estimators appear as generalized maximum likelihood estimators and ensure good asymptotic properties~\cite{pascal2008performance, OTKP12a, mahot2013asymptotic}.
Nevertheless, $M$-estimators do not account for the low-rank structure.
A natural solution to this issue is to directly derive an estimator as the minimizer of a robust cost function under a low-rank structure constraint. 
This approach has been proposed in~\cite[Sec. V.A.]{SBP16}, where a majorization-minimization algorithms is proposed to minimize Tyler's cost function according to this structure.
However, the tractability of this estimator is an open question at low sample support (\emph{cf.} assumption 2 in~\cite{SBP16}). 
Notably, the majorization-minimization algorithm can present convergence issues in some practical case where $n$ is close to or smaller than $p$. 

To address this issue, we propose to use the Riemannian optimization framework~\cite{AMS08}: the proposed geometry for the the quotient $(\textup{St}_{p,k}\times\mathcal{H}^{++}_k)/\mathcal{U}_k$ indeed offers the possibility to apply a large panel of generic first and second order optimization algorithms on manifolds, such as gradient descent, conjugate gradient, BFGS, trust region, Newton, \emph{etc.} (\emph{cf.}~\cite{AMS08} for details).
More specifically for robust covariance matrix estimation, we propose an estimator formulated as the minimizer of a counterpart of Tyler's cost function defined directly on $(\textup{St}_{p,k}\times\mathcal{H}^{++}_k)/\mathcal{U}_k$.
We then focus on two algorithms for solving the introduced problem: one based on Riemannian gradient descent (first order method), the other based on Riemannian trust region (second order method).
In terms of estimation accuracy, our numerical experiments show that the Riemannian trust region based algorithm is similar to~\cite[algorithm 5]{SBP16}.
Interestingly, these experiments also show that the Riemannian gradient descent based method can still reach good performance when the other methods diverge at insufficient sample support.  

\subsection{Statistical performance analysis in low-rank elliptical models}

Cram\'er-Rao lower bounds are ubiquitous tools to characterize the optimum performances in terms of mean squared error that can be achieved for a given parametric estimation problem~\cite{K93}. 
In the context of elliptical distributions, Cram\'er-Rao lower bounds can be obtained using the general results of~\cite{BA13}, and have been studied for covariance/shape estimation in~\cite{GG13, PR10}. 
However, the low-rank models involve constraints and ambiguities on the parameters space, which does not allow for simple/practical derivations, even using the so-called constrained Cram\'er-Rao lower bounds~\cite{GH90, M93, SN98}.
Additionally, the classical inequality applies on the mean squared error (Euclidean metric), while this criterion may not be the most appropriate for characterizing the performance when parameters are living in a manifold.

To overcome these issues, we consider the framework of intrinsic Cram\'er-Rao lower bounds from~\cite{S05, B13, BGRB18}.
For covariance matrix estimation in low-rank elliptical models, two performance criteria are considered:
the proposed divergence on the quotient $(\textup{St}_{p,k}\times\mathcal{H}^{++}_k)/\mathcal{U}_k$ (for total error measurement), and the Riemmanian distance on the Grassmann manifold $\mathcal{G}_{p,k}$~\cite{EAS98} (for principal subspace estimation error measurement).
We derive lower bounds for both error measures and observe thanks to numerical experiments that they can be reached by the proposed algorithms. 
These contributions therefore generalize the ones of~\cite{S05} on low-rank Gaussian models to wider classes of distributions and performance measures.

\section{Background}
\label{sec:model}
\subsection{Complex elliptically symmetric distributions and robust covariance estimation}

Complex elliptically symmetric distributions~\cite{KTYT90} represent a large family of multivariate distributions that encompasses, for example, Gaussian, $K$-, Student $t$-, and Weibull distributions.
A detailed review on the topic can be found in~\cite{OTKP12a}.
The probability density function (pdf) associated with the random variable $\VEC{x}\in\mathbb{C}^p$ following a zero-mean complex elliptically symmetric distribution is, up to a normalization factor,
\begin{equation}
	f_g^{++}(\VEC{x}|\MAT{R}) = \det(\MAT{R})^{-1} g(\VEC{x}^H \MAT{R}^{-1}\VEC{x}),
\label{eq:model:ces_pdf_HPD}
\end{equation}
where $\det$ denotes the determinant operator, $\MAT{R}\in\mathcal{H}^{++}_p$ is the covariance matrix and $g:\mathbb{R}^{+}\to\mathbb{R}^{+}$ is the so-called density generator of the distribution.

The negative log-likelihood function associated with $n$ independent and identically distributed samples $\{\VEC{x}_i\}$ of the random variable $\VEC{x}$ is
\begin{equation*}
	L_g^{++}(\MAT{R}) = n\log\det(\MAT{R}) - \sum_{i=1}^n \log(g(\VEC{x}_i^H\MAT{R}^{-1} \VEC{x}_i)).
\end{equation*}
Given the density generator $g$ and $n$ observations $\{\VEC{x}_i\}$, an estimator $\MAT{\widehat{R}}$ of the true covariance matrix $\MAT{R}$ can be obtained by solving the optimization problem
\begin{equation*}
	\MAT{\widehat{R}} \, = \, \argmin_{\MAT{R}} \quad L_g^{++}(\MAT{R}).
\end{equation*}
Unfortunately, the true density generator $g$ is often unknown in practice.
To overcome this issue, a solution provided by the robust estimation theory is to compute an $M$-estimator~\cite{MY76}.
A popular choice is Tyler's $M$-estimator~\cite{T87, pascal2008covariance}, which is motivated by its ``distribution-free'' properties among the whole familly of CES, its good asymptotic performance~\cite{ pascal2008covariance}, and robustness properties.
Given $\{\VEC{x}_i\}$, the corresponding cost function to be minimized corresponds to $g(t)=1/t$ and is defined as
\begin{equation}
	L^{++}_{\textup{T}}(\MAT{R}) = p\sum_{i=1}^n \log(\VEC{x}_i^H\MAT{R}^{-1}\VEC{x}_i)	+ n \log\det(\MAT{R}).
\label{eq:model:tyler_HPD}
\end{equation}
On $\mathcal{H}^{++}_p$, this cost function is efficiently minimized with a fixed-point algorithm~\cite{pascal2008covariance}.
Additional assumptions on the structure of the covariance $\MAT{R}$ can also be made; see \emph{e.g.},~\cite{SBP16} for various possibilities.
In this work, we are interested in the low-rank covariance structure, which is for instance treated in~\cite[section V.A]{SBP16} and~\cite{S05}.

\subsection{Low-rank covariance model and parameter space}
\label{subsec:model:cov_hlr}

The low-rank covariance model (also known as spiked model~\cite{J01} or factor model~\cite{R11}) refers to the structure%
\footnote{
	One might be interested in the more general model $\MAT{R}=\MAT{R}_0+\MAT{H}$, where the identity $\MAT{I}_p$ is replaced by any (known) $\MAT{R}_0\in\mathcal{H}^{++}_p$, as done in~\cite{S05}.
	It is equivalent to our model as it suffices to whiten the random variable $\VEC{x}$ with $\MAT{\Sigma}_0^{\nicefrac{-1}2}$ in order to obtain~\eqref{eq:model:spiked_model}.
	Furthermore, as done in many works, we assume the rank $k$ to be known (\emph{e.g.}, from prior physical considerations~\cite{GS07}) or pre-estimated (\emph{e.g.}, from model order selection techniques~\cite{SS04}).
}
\begin{equation}
\MAT{R} = \MAT{I}_p + \MAT{H},
\label{eq:model:spiked_model}
\end{equation}
where $\MAT{I}_p$ denotes the $p$-dimensional identity matrix and $\MAT{H}$ is a $p\times p$ Hermitian positive semi-definite matrix of rank~$k$.
This model is directly related to principal component analysis and subspace recovery~\cite{TB99}.
Even though a Majorization-Minimization algorithm is proposed in~\cite{SBP16} to treat this particular problem, the tractability of the resulting estimator is an open question for $n<p$ (\emph{cf.} \cite[assumption 2]{SBP16}), and convergence issues are observed in some practical cases.

The parameter $\MAT{H}$ in~\eqref{eq:model:spiked_model} lives in the manifold $\mathcal{H}^{+}_{p,k}$ of $p\times p$ Hermitian positive semi-definite matrices of rank $k$.
As explained in the introduction, several geometries have been proposed for this manifold.
In this work, we consider the geometry resulting from the decomposition
\begin{equation}
	\MAT{H}=\MAT{U}\MAT{\Sigma}\MAT{U}^H, ~\textup{with}~ (\MAT{U},\MAT{\Sigma})\in\overbar{\mathcal{M}}_{p,k}=(\textup{St}_{p,k}\times\mathcal{H}^{++}_k),
\label{eq:model:HLR_decomp}
\end{equation}
which is directly related to the singular value decomposition of $\MAT{H}$.
This parametrization is particularly interesting when it comes to subspace estimation as the latter is simply obtained from the component $\MAT{U}$.

Let $\overbar{\varphi}:\overbar{\mathcal{M}}_{p,k}\to\mathcal{H}^+_{p,k}$ be the smooth mapping defined, for $(\MAT{U},\MAT{\Sigma})\in\overbar{\mathcal{M}}_{p,k}$, as
\begin{equation}
	\overbar{\varphi}(\MAT{U},\MAT{\Sigma})=\MAT{U}\MAT{\Sigma}\MAT{U}^H.
\label{eq:model:mapping_product2HLR}
\end{equation}
Since every $\MAT{H}\in\mathcal{H}^{+}_{p,k}$ admits a decomposition of the form~\eqref{eq:model:HLR_decomp}, the mapping $\overbar{\varphi}$ is surjective.
However, it is not injective as the considered decomposition is not unique: given any $\MAT{O}\in\mathcal{U}_k$, one has $\MAT{H}=\overbar{\varphi}(\MAT{U},\MAT{\Sigma})=\overbar{\varphi}(\MAT{U}\MAT{O},\MAT{O}^H\MAT{\Sigma}\MAT{O})$.
As done in~\cite{BS09,MBS11}, to account for the action of the unitary matrices, we define the quotient manifold
\begin{equation}
	\mathcal{M}_{p,k} = \{ \pi(\MAT{U},\MAT{\Sigma}): \, (\MAT{U},\MAT{\Sigma})\in\overbar{\mathcal{M}}_{p,k} \},
\label{eq:model:quotient}
\end{equation}
where the equivalence class $\pi(\MAT{U},\MAT{\Sigma})$ is
\begin{equation}
	\pi(\MAT{U},\MAT{\Sigma}) = \{ (\MAT{U}\MAT{O},\MAT{O}^H\MAT{\Sigma}\MAT{O}) : \, \MAT{O}\in\mathcal{U}_k \}.
\label{eq:model:eq_class}
\end{equation}

As shown in~\cite{BS09,MBS11}, it follows that the function $\varphi$ on $\mathcal{M}_{p,k}$ induced by $\overbar{\varphi}$ on $\overbar{\mathcal{M}}_{p,k}$, \emph{i.e.}, such that $\overbar{\varphi}=\varphi\circ\pi$, is an isomorphism from $\mathcal{M}_{p,k}$ onto $\mathcal{H}^+_{p,k}$.
Thus, the geometry of $\mathcal{M}_{p,k}$ can be exploited to treat problems defined on $\mathcal{H}^+_{p,k}$.
In particular, the pdf on $\mathcal{M}_{p,k}$ of a random variable $\VEC{x}$ following a zero-mean complex elliptically symmetric distribution with covariance matrix admitting structure~\eqref{eq:model:HLR_decomp} is, for all $\theta=\pi(\MAT{U},\MAT{\Sigma})\in\mathcal{M}_{p,k}$,
\begin{equation}
	f_g(\VEC{x}|\theta) = f_g^{++}(\VEC{x}|\MAT{I}_p+\varphi(\theta)),
\label{eq:model:ces_pdf_HLR}
\end{equation}
where $f_g^{++}$ is defined in~\eqref{eq:model:ces_pdf_HPD}.
Similarly, the cost function on $\mathcal{M}_{p,k}$ of the Tyler's $M$-estimator is defined, for all $\theta=\pi(\MAT{U},\MAT{\Sigma})\in\mathcal{M}_{p,k}$, as
\begin{equation}
	L_{\textup{T}}(\theta) = L^{++}_{\textup{T}}(\MAT{I}_p + \varphi(\theta)),
\label{eq:model:tyler_HLR}
\end{equation}
where $L^{++}_{\textup{T}}$ is defined in~\eqref{eq:model:tyler_HPD}.

\section{Riemannian geometry of Hermitian positive semi-definite matrices of fixed rank}
\label{sec:rg_hlr}

\begin{figure}
	\centering
	\begin{tikzpicture}[scale=3.5]
\draw [fill=gray!0] (1,0) arc (30:150:1.155)
      plot [smooth, domain=pi:2*pi] ({cos(\x r)},{0.2*sin(\x r)});
\draw [dotted] plot [smooth, domain=0:pi] ({cos(\x r)},{0.2*sin(\x r)});
\draw (1.08,0.15) node {$\overbar{\mathcal{M}}_{p,k}$};

\draw (-0.3,0.535) to[bend right=30] (-0.8,-0.125);
\draw (0,0.58) to[bend right=20] (-0.4,-0.18);
\draw (0.25,0.55) to[bend right=15] (-0.05,-0.2);
\draw (0.58,0.425) to[bend right=10] (0.34,-0.19);

\coordinate (x) at (-0.56,0.37);

\draw [fill=gray!20,opacity=0.4] (x)++(0.55,0.07) -- ++(-0.5,-0.4) -- ++(-0.5,+0.25) -- ++(+0.5,+0.4) -- cycle;
\draw [opacity=0.8] (-0.69,0.73) node {$T_{\overbar{\theta}}\overbar{\mathcal{M}}_{p,k}$};

\draw[draw=black!70,line width=1.2pt] (-0.312,0.592) -- (-0.782,0.175);
\draw [opacity=0.7] (-0.3,0.67) node {$\mathcal{V}_{\overbar{\theta}}$};
\draw[draw=black!70,line width=1.2pt] (-0.78,0.48) -- (-0.27,0.23);
\draw [opacity=0.7] (-0.85,0.53) node {$\mathcal{H}_{\overbar{\theta}}$};

\draw (x) node {$\bullet$};
\draw (x)++(-0.02,+0.04) node[above] {\small$\overbar{\theta}$};

\draw (-0.77,-0.195) node[rotate=-10] {\footnotesize$\pi^{-1}(\pi(\overbar{\theta}))$};

\draw[>=stealth,->] (x) -- ++(0.21,-0.1) node [below,midway,sloped] {\small$P^{\mathcal{H}}_{\overbar{\theta}}(\overbar{\xi})$};
\draw[>=stealth,->] (x) -- ++(0.3,0.0);
\draw (-0.28,0.45) node {\small$\overbar{\xi}$};
\draw[dashed,opacity=0.8] (-0.35,0.27) -- (-0.26,0.37);

\draw[>=stealth,->,line width=1.5pt] (-0.1,-0.28) -- (-0.1,-0.45);
\draw (0,-0.35)  node {$\pi$};

\draw (-1.1,-0.72) to[bend left=20] (1,-0.72);
\draw (1.05,-0.62) node {$\mathcal{M}_{p,k}$};
\draw (-0.75,-0.7) node {$\bullet$};
\draw (-0.75,-0.78) node {\small$\theta=\pi(\overbar{\theta})$};

\draw (-0.35,-0.7) node {$\bullet$};
\draw (0,-0.7) node {$\bullet$};
\draw (0.4,-0.7) node {$\bullet$};

\end{tikzpicture}
	\caption{
		Illustration of the quotient manifold $\mathcal{M}_{p,k}$ of the manifold $\overbar{\mathcal{M}}_{p,k}$.
		The tangent space $T_{\overbar{\theta}}\overbar{\mathcal{M}}_{p,k}$ can be decomposed into two complementary subspaces: the vertical space $\mathcal{V}_{\overbar{\theta}}=T_{\overbar{\theta}}\pi^{-1}(\pi(\overbar{\theta}))$ and the horizontal space $\mathcal{H}_{\overbar{\theta}}$, which provides proper representatives of tangent vectors in $T_{\theta}\mathcal{M}_{p,k}$ at $\theta=\pi(\overbar{\theta})$.
		The orthogonal projection map $P^{\mathcal{H}}_{\overbar{\theta}}$ allows to project $\overbar{\xi}\in T_{\overbar{\theta}}\overbar{\mathcal{M}}_{\theta}$ onto $\mathcal{H}_{\overbar{x}}$.
		Both $\mathcal{H}_{\overbar{\theta}}$ and $P^{\mathcal{H}}_{\overbar{\theta}}$ are defined in proposition~\ref{prop:rg_hlr:horizontal}.
	}
\label{fig:quotient}
\end{figure}
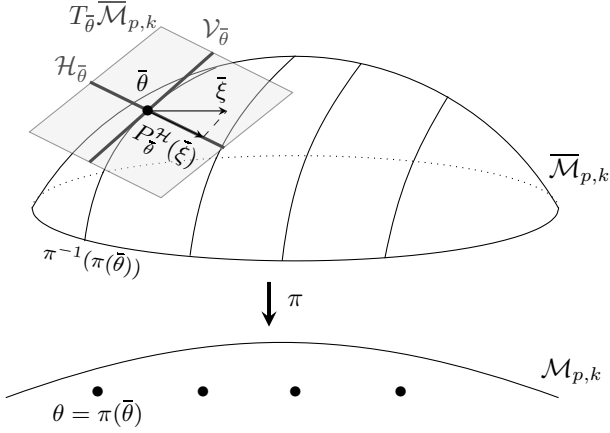

To describe the geometry of the quotient $\mathcal{M}_{p,k}$, we exploit the submersion $\pi:\overbar{\mathcal{M}}_{p,k}\to\mathcal{M}_{p,k}$ defined in~\eqref{eq:model:eq_class}.
This allows to work with representatives of the geometrical objects of the quotient in $\overbar{\mathcal{M}}_{p,k}$.
In particular, $\theta\in\mathcal{M}_{p,k}$ is represented by any $\overbar{\theta}=(\MAT{U},\MAT{\Sigma})\in\overbar{\mathcal{M}}_{p,k}$ such that $\theta=\pi(\overbar{\theta})$.
The tangent space $T_{\theta}\mathcal{M}_{p,k}$ at $\theta=\pi(\overbar{\theta})$ in $\mathcal{M}_{p,k}$ is represented by a well chosen subspace of the tangent space $T_{\overbar{\theta}}\overbar{\mathcal{M}}_{p,k}$ at $\overbar{\theta}$ in $\overbar{\mathcal{M}}_{p,k}$.
Moreover, a Riemannian metric on $\mathcal{M}_{p,k}$ can be defined through a metric on $\overbar{\mathcal{M}}_{p,k}$ that is invariant along the equivalence classes~\eqref{eq:model:eq_class}.
An illustration of the quotient $\mathcal{M}_{p,k}$ is provided in figure~\ref{fig:quotient}.

In the following, $\overbar{\theta}=(\MAT{U},\MAT{\Sigma})$, $\overbar{\xi}=(\MAT{\xi}_{\MAT{U}},\MAT{\xi}_{\MAT{\Sigma}})$, $\overbar{\eta}=(\MAT{\eta}_{\MAT{U}},\MAT{\eta}_{\MAT{\Sigma}})$ and $\overbar{Z}=(\MAT{Z}_{\MAT{U}},\MAT{Z}_{\MAT{\Sigma}})$.
First recall that
\begin{equation}
	T_{\overbar{\theta}}\overbar{\mathcal{M}}_{p,k} =
	\{ \overbar{\xi}\in\mathbb{C}^{p\times k}\times\mathcal{H}_k : \, \MAT{U}^H\MAT{\xi}_{\MAT{U}} + \MAT{\xi}_{\MAT{U}}^H\MAT{U} = \MAT{0} \}.
\label{eq:rg_hlr:tangent_product}
\end{equation}
We equip $\overbar{\mathcal{M}}_{p,k}$ with the Riemannian metric of definition~\ref{def:rg_hlr:metric}.
The part of this metric that concerns $\MAT{U}$ is the so-called canonical metric on Stiefel~\cite{EAS98}%
\footnote{
	This metric is advantageous as compared to the Euclidean metric because resulting geodesics admit simpler formulas~\cite{EAS98}.
},
which is obtained by treating $\textup{St}_{p,k}$ as the quotient $\mathcal{U}_p/\mathcal{U}_{p-k}$.
The one that concerns $\MAT{\Sigma}$ corresponds to a class of affine invariant metrics on $\mathcal{H}^{++}_k$ that are of interest when dealing with elliptical distributions as they are related to the Fisher information metric~\cite{BGRB18}%
\footnote{
	For example, the Fisher information metric on $\mathcal{H}^{++}_k$ for the Gaussian distribution is obtained with $\alpha=1$ and $\beta=0$.
}.
\begin{definition}[Riemannian metric]
\label{def:rg_hlr:metric}
	We define the Riemannian metric $\langle\cdot,\cdot\rangle_{\cdot}$ on $\overbar{\mathcal{M}}_{p,k}$ by
	\begin{multline}
		\langle \overbar{\xi},\overbar{\eta} \rangle_{\overbar{\theta}} = 
		\reel(\tr(\MAT{\xi}_{\MAT{U}}^H(\MAT{I}_p - \frac12\MAT{U}\MAT{U}^H)\MAT{\eta}_{\MAT{U}}))
		\\
		+ \alpha\tr(\MAT{\Sigma}^{-1}\MAT{\xi}_{\MAT{\Sigma}}\MAT{\Sigma}^{-1}\MAT{\eta}_{\MAT{\Sigma}})+\beta\tr(\MAT{\Sigma}^{-1}\MAT{\xi}_{\MAT{\Sigma}})\tr(\MAT{\Sigma}^{-1}\MAT{\eta}_{\MAT{\Sigma}}),
	\label{eq:rg_hlr:metric}
	\end{multline}
	where $\alpha>0$ and $\beta>-\frac{\alpha}{k}$.
\end{definition}
It is readily checked that the metric~\eqref{eq:rg_hlr:metric} is invariant along the equivalence classes~\eqref{eq:model:eq_class}, \emph{i.e.}, for all $\MAT{O}\in\mathcal{U}_k$
\begin{equation*}
	\langle \overbar{\xi}, \overbar{\eta} \rangle_{\overbar{\theta}} = \langle \phi_{\MAT{O}}(\overbar{\xi}), \phi_{\MAT{O}}(\overbar{\eta}) \rangle_{\phi_{\MAT{O}}(\overbar{\theta})},
\end{equation*}
where $\phi_{\MAT{O}}(\overbar{Z})=(\MAT{Z}_{\MAT{U}}\MAT{O},\MAT{O}^H\MAT{Z}_{\MAT{\Sigma}}\MAT{O})$.
Thus, metric~\eqref{eq:rg_hlr:metric} induces a Riemannian metric on the quotient $\mathcal{M}_{p,k}$.
Furthermore, the orthogonal projection map according to~\eqref{eq:rg_hlr:metric} from $\mathbb{C}^{p\times k}\times\mathbb{C}^{k\times k}$ onto $T_{\overbar{\theta}}\overbar{\mathcal{M}}_{p,k}$~is
\begin{equation}
	P_{\overbar{\theta}}(\overbar{Z}) = (\MAT{Z}_{\MAT{U}}-\MAT{U}\herm(\MAT{U}^H\MAT{Z}_{\MAT{U}}), \herm(\MAT{Z}_{\MAT{\Sigma}})),
\label{eq:rg_hlr:proj_product}
\end{equation}
where $\herm$ returns the Hermitian part of its argument.

The tangent space $T_{\overbar{\theta}}\overbar{\mathcal{M}}_{p,k}$ can be decomposed into two complementary spaces: the vertical and horizontal spaces $\mathcal{V}_{\overbar{\theta}}$ and $\mathcal{H}_{\overbar{\theta}}$~\cite{AMS08}.
The vertical space is the tangent space $T_{\overbar{\theta}}\pi^{-1}(\pi(\overbar{\theta}))$ to the equivalence class $\pi^{-1}(\pi(\overbar{\theta}))$ at $\overbar{\theta}$, which, as shown in~\cite{BS09,MBS11}, is given~by
\begin{equation*}
	\mathcal{V}_{\overbar{\theta}} = \{ (\MAT{U}\MAT{\Omega}, \MAT{\Sigma}\MAT{\Omega} - \MAT{\Omega}\MAT{\Sigma}): \, \MAT{\Omega}\in\mathcal{H}^{\perp}_k\},
\end{equation*}
where $\mathcal{H}^{\perp}_k$ denotes the space of skew-Hermitian matrices.
$\mathcal{H}_{\overbar{\theta}}$, which provides proper representatives for the elements of $T_{\theta}\mathcal{M}_{p,k}$%
\footnote{
	Given $\theta=\pi(\overbar{\theta})\in\mathcal{M}_{p,k}$, the tangent vector $\xi\in T_{\theta}\mathcal{M}_{p,k}$ is represented by the only $\overbar{\xi}\in\mathcal{H}_{\overbar{\theta}}$ such that $\xi=\D\pi(\overbar{\theta})[\overbar{\xi}]$.
}
and turns $\pi$ into a Riemannian submersion, is then defined as the orthogonal complement to $\mathcal{V}_{\overbar{\theta}}$ according to metric~\eqref{eq:rg_hlr:metric}.
The horizontal space along with the orthogonal projection map from $T_{\overbar{\theta}}\overbar{\mathcal{M}}_{p,k}$ onto $\mathcal{H}_{\overbar{\theta}}$ are given in proposition~\ref{prop:rg_hlr:horizontal}.
\begin{proposition}
\label{prop:rg_hlr:horizontal}
	The horizontal space $\mathcal{H}_{\overbar{\theta}}$ at $\overbar{\theta}\in\overbar{\mathcal{M}}_{p,k}$ is
	\begin{equation*}
		\mathcal{H}_{\overbar{\theta}} = \{\overbar{\xi}\in T_{\overbar{\theta}}\overbar{\mathcal{M}}_{p,k}: \, \MAT{U}^H\MAT{\xi}_{\MAT{U}} = 2\alpha(\MAT{\Sigma}^{-1}\MAT{\xi}_{\MAT{\Sigma}} - \MAT{\xi}_{\MAT{\Sigma}}\MAT{\Sigma}^{-1}) \}.
	\end{equation*}
	The orthogonal projection map $P^{\mathcal{H}}_{\overbar{\theta}}$ according to~\eqref{eq:rg_hlr:metric} from $T_{\overbar{\theta}}\overbar{\mathcal{M}}_{p,k}$ onto $\mathcal{H}_{\overbar{\theta}}$ is given by
	\begin{equation*}
		P^{\mathcal{H}}_{\overbar{\theta}}(\overbar{\xi}) = (\MAT{\xi}_{\MAT{U}} - \MAT{U}\MAT{\Omega}, \MAT{\xi}_{\MAT{\Sigma}} + \MAT{\Omega}\MAT{\Sigma} - \MAT{\Sigma}\MAT{\Omega}),
	\end{equation*}
	where $\MAT{\Omega}\in\mathcal{H}^{\perp}_k$ is the unique solution to
	\begin{multline*}
		(1-4\alpha)\MAT{\Omega} + 2\alpha(\MAT{\Sigma}^{-1}\MAT{\Omega}\MAT{\Sigma} + \MAT{\Sigma}\MAT{\Omega}\MAT{\Sigma}^{-1}) = 
		\\
		\MAT{U}^H\MAT{\xi}_{\MAT{U}} + 2\alpha(\MAT{\xi}_{\MAT{\Sigma}}\MAT{\Sigma}^{-1} + \MAT{\Sigma}^{-1}\MAT{\xi}_{\MAT{\Sigma}}).
	\end{multline*}
\end{proposition}
\begin{IEEEproof}
	By definition, $\overbar{\xi}\in\mathcal{H}_{\overbar{\theta}}$ if and only if, for all $\MAT{\Omega}\in\mathcal{H}^{\perp}_k$, $\langle \overbar{\xi}, (\MAT{U}\MAT{\Omega}, \MAT{\Sigma}\MAT{\Omega} - \MAT{\Omega}\MAT{\Sigma}) \rangle_{\overbar{\theta}} = 0$.
	From~\eqref{eq:rg_hlr:metric}, basic calculations yield $\tr((\MAT{\xi}_{\MAT{U}}^H\MAT{U} + 2\alpha(\MAT{\Sigma}^{-1}\MAT{\xi}_{\MAT{\Sigma}} - \MAT{\xi}_{\MAT{\Sigma}}\MAT{\Sigma}^{-1}))\MAT{\Omega}) = 0$.
	This is true for all $\MAT{\Omega}\in\mathcal{H}^{\perp}_k$ if and only if $\MAT{\xi}_{\MAT{U}}^H\MAT{U} + 2\alpha(\MAT{\Sigma}^{-1}\MAT{\xi}_{\MAT{\Sigma}} - \MAT{\xi}_{\MAT{\Sigma}}\MAT{\Sigma}^{-1})$ is Hermitian.
	This translates into $\MAT{U}^H\MAT{\xi}_{\MAT{U}} - \MAT{\xi}_{\MAT{U}}^H\MAT{U} = 4\alpha(\MAT{\Sigma}^{-1}\MAT{\xi}_{\MAT{\Sigma}} - \MAT{\xi}_{\MAT{\Sigma}}\MAT{\Sigma}^{-1})$.
	From~\eqref{eq:rg_hlr:tangent_product}, we have $\MAT{U}^H\MAT{\xi_U} + \MAT{\xi_U}^H\MAT{U} = \MAT{0}$, leading to the result.

	Regarding $P^{\mathcal{H}}$, it has the proposed form by definition.
	The matrix $\MAT{\Omega}\in\mathcal{H}^{\perp}_k$ must be chosen in order to have $P^{\mathcal{H}}_{\overbar{\theta}}(\overbar{\xi})\in\mathcal{H}_{\overbar{\theta}}$.
	Basic calculations yield the proposed equation.
	It remains to show that the solution exists and is unique.
	This equation can be vectorized as
	\begin{multline*}
		((1-4\alpha)\MAT{I}_{k^2} + 2\alpha(\MAT{\Sigma}^{-T}\otimes\MAT{\Sigma} + \MAT{\Sigma}^T\otimes\MAT{\Sigma}^{-1}))\vect(\MAT{\Omega}) =
		\\
		\vect(\MAT{U}^H\MAT{\xi}_{\MAT{U}} + 2\alpha(\MAT{\xi}_{\MAT{\Sigma}}\MAT{\Sigma}^{-1} + \MAT{\Sigma}^{-1}\MAT{\xi}_{\MAT{\Sigma}})).
	\end{multline*}
	Showing that $(1-4\alpha)\MAT{I}_{p^2} + 2\alpha(\MAT{\Sigma}^{-T}\otimes\MAT{\Sigma} + \MAT{\Sigma}^T\otimes\MAT{\Sigma}^{-1})$ is positive definite is enough to conclude.
	In order to do so, consider the eigenvalue decomposition $\MAT{\Sigma}=\MAT{V}\MAT{\Lambda}\MAT{V}^H$.
	We have
	\begin{multline*}
		((1-4\alpha)\MAT{I}_{k^2} + 2\alpha(\MAT{\Sigma}^{-1}\otimes\MAT{\Sigma} + \MAT{\Sigma}\otimes\MAT{\Sigma}^{-1})) =
		\\
		(\MAT{\overbar{V}}\otimes\MAT{V})((1-4\alpha)\MAT{I}_{k^2} + 2\alpha(\MAT{\Lambda}^{-1}\otimes\MAT{\Lambda} + \MAT{\Lambda}\otimes\MAT{\Lambda}^{-1}))(\MAT{\overbar{V}}\otimes\MAT{V})^H,
	\end{multline*}
	where $\MAT{\overbar{V}}$ is the conjugate of $\MAT{V}$.
	As $\MAT{V}$ is unitary, $\MAT{\overbar{V}}$ and $\MAT{\overbar{V}}\otimes\MAT{V}$ are also unitary.
	$((1-4\alpha)\MAT{I}_{p^2} + 2\alpha(\MAT{\Lambda}^{-1}\otimes\MAT{\Lambda} + \MAT{\Lambda}\otimes\MAT{\Lambda}^{-1}))$ is diagonal and its elements are $1 - 4\alpha + 2\alpha(\frac{\lambda_i}{\lambda_j} + \frac{\lambda_j}{\lambda_i})$, where $\lambda_i$ is the $i^{\textup{th}}$ diagonal element of $\MAT{\Lambda}$.
	The function $h(x)=x+\frac1x$, defined for $x>0$, admits $2$ as a global minimum for $x=1$, showing that $1 - 4\alpha + 2\alpha(\frac{\lambda_i}{\lambda_j} + \frac{\lambda_j}{\lambda_i})\geq 1>0$.
	This completes the proof.
\end{IEEEproof}

\vspace*{.5cm}
The Levi-Civita connection on $\mathcal{M}_{p,k}$ associated with the metric induced by~\eqref{eq:rg_hlr:metric}, which generalizes the concept of directional derivative of vector fields on a manifold%
\footnote{
	A vector field is an operator which assigns a tangent vector to every point of a manifold.
	An example of a vector field is the gradient of an objective function.
},
is given in proposition~\ref{prop:rg_hlr:LC}.
This object is crucial when it comes to defining geodesics and the Riemannian Hessian of an objective function on $\mathcal{M}_{p,k}$.
\begin{proposition}
\label{prop:rg_hlr:LC}
	Let $\theta=\pi(\overbar{\theta})\in\mathcal{M}_{p,k}$, $\xi=\D\pi(\overbar{\theta})[\overbar{\xi}]\in T_{\theta}\mathcal{M}_{p,k}$ and the vector field $\eta=\D\pi(\overbar{\theta})[\overbar{\eta}]$ evaluated at $\theta$, where $\overbar{\xi}$, $\overbar{\eta}\in\mathcal{H}_{\overbar{\theta}}$.
	The representative $\overbar{\nabla_{\xi}\, \eta}$ in $\mathcal{H}_{\overbar{\theta}}$ of the Levi-Civita connection $\nabla_{\xi}\, \eta$ on $\mathcal{M}_{p,k}$ is
	\begin{equation*}
		\overbar{\nabla_{\xi}\, \eta} = P^{\mathcal{H}}_{\overbar{\theta}}(\overbar{\nabla}_{\overbar{\xi}}\, \overbar{\eta}),
	\end{equation*}
	where $\overbar{\nabla}_{\overbar{\xi}}\, \overbar{\eta}$ is the Levi-Civita connection on $\overbar{\mathcal{M}}_{p,k}$, given by
	\begin{multline*}
		\overbar{\nabla}_{\overbar{\xi}}\, \overbar{\eta} = P_{\overbar{\theta}}(\D\overbar{\eta}[\overbar{\xi}])
		+ ((\MAT{I}_p-\MAT{U}\MAT{U}^H)\herm(\MAT{\eta}_{\MAT{U}}\MAT{\xi}_{\MAT{U}}^H)\MAT{U}, \\ -\herm(\MAT{\eta}_{\MAT{\Sigma}}\MAT{\Sigma}^{-1}\MAT{\xi}_{\MAT{\Sigma}}))
	\end{multline*}
\end{proposition}
\begin{IEEEproof}
	Let $\overbar{g}_{\overbar{\theta}}(\overbar{\xi},\overbar{\eta}) = \langle\overbar{\xi},\overbar{\eta}\rangle_{\overbar{\theta}}$.
	The Koszul formula~\cite{AMS08}, which characterizes the Levi-Civita connection, is in our case
	\begin{multline*}
		2\overbar{g}_{\overbar{\theta}}(\overbar{\nabla}_{\overbar{\xi}}\,\overbar{\eta},\overbar{\nu}) - 2\overbar{g}_{\overbar{\theta}}(\D\overbar{\eta}[\overbar{\xi}],\overbar{\nu}) =
		\\
		+ \D\overbar{g}_{\overbar{\theta}}[\overbar{\xi}](\overbar{\eta},\overbar{\nu})
		+ \D\overbar{g}_{\overbar{\theta}}[\overbar{\eta}](\overbar{\xi},\overbar{\nu})
		- \D\overbar{g}_{\overbar{\theta}}[\overbar{\nu}](\overbar{\xi},\overbar{\eta}).
	\end{multline*}
	To obtain the three terms on the right side of this equation, we have to derive the metric $\overbar{g}_{\overbar{\theta}}$ with respect to $\overbar{\theta}$.
	One can check that
	\begin{multline*}
		\D\overbar{g}_{\overbar{\theta}}[\overbar{\nu}](\overbar{\xi},\overbar{\eta}) =
		-\reel(\tr(\MAT{\xi}_{\MAT{U}}^H\herm(\MAT{U}\MAT{\nu}_{\MAT{U}}^H)\MAT{\eta}_{\MAT{U}}))
		\\
		-\beta\tr(\MAT{\Sigma}^{-1}\MAT{\xi}_{\MAT{\Sigma}}\MAT{\Sigma}^{-1}\MAT{\nu}_{\MAT{\Sigma}})\tr(\MAT{\Sigma}^{-1}\MAT{\eta}_{\MAT{\Sigma}})
		\\
		-\beta\tr(\MAT{\Sigma}^{-1}\MAT{\xi}_{\MAT{\Sigma}})\tr(\MAT{\Sigma}^{-1}\MAT{\eta}_{\MAT{\Sigma}}\MAT{\Sigma}^{-1}\MAT{\nu}_{\MAT{\Sigma}})
		\\
		-2\alpha\tr(\MAT{\Sigma}^{-1}\MAT{\xi}_{\MAT{\Sigma}}\MAT{\Sigma}^{-1}\MAT{\eta}_{\MAT{\Sigma}}\MAT{\Sigma}^{-1}\MAT{\nu}_{\MAT{\Sigma}}).
	\end{multline*}
	It follows that the right side of the Koszul formula is
	\begin{multline*}
		\D\overbar{g}_{\overbar{\theta}}[\overbar{\xi}](\overbar{\eta},\overbar{\nu})
		+ \D\overbar{g}_{\overbar{\theta}}[\overbar{\eta}](\overbar{\xi},\overbar{\nu})
		- \D\overbar{g}_{\overbar{\theta}}[\overbar{\nu}](\overbar{\xi},\overbar{\eta})
		=
		\\
		\tr(\MAT{\nu}_{\MAT{U}}^H(2\herm(\MAT{\eta}_{\MAT{U}}\MAT{\xi}_{\MAT{U}}^H)\MAT{U}-\MAT{U}\herm(\MAT{\eta}_{\MAT{U}}^H\MAT{\xi}_{\MAT{U}})))
		\\
		-2\alpha\tr(\MAT{\Sigma}^{-1}\MAT{\xi}_{\MAT{\Sigma}}\MAT{\Sigma}^{-1}\MAT{\eta}_{\MAT{\Sigma}}\MAT{\Sigma}^{-1}\MAT{\nu}_{\MAT{\Sigma}})
		\\
		-2\beta\tr(\MAT{\Sigma}^{-1}\MAT{\xi}_{\MAT{\Sigma}}\MAT{\Sigma}^{-1}\MAT{\eta}_{\MAT{\Sigma}})\tr(\MAT{\Sigma}^{-1}\MAT{\nu}_{\MAT{\Sigma}}).
	\end{multline*}
	Moreover,
	\begin{equation*}
		\tr(\MAT{\nu}_{\MAT{U}}^H\MAT{\widetilde{Z}}_{\MAT{U}}) = \tr(\MAT{\nu}_{\MAT{U}}^H(\MAT{I}_p -\frac12\MAT{U}\MAT{U}^H)(\MAT{I}_p+\MAT{U}\MAT{U}^H)\MAT{\widetilde{Z}}_{\MAT{U}}).
	\end{equation*}
	It follows that
	\begin{equation*}
		\D\overbar{g}_{\overbar{\theta}}[\overbar{\xi}](\overbar{\eta},\overbar{\nu})
		+ \D\overbar{g}_{\overbar{\theta}}[\overbar{\eta}](\overbar{\xi},\overbar{\nu})
		- \D\overbar{g}_{\overbar{\theta}}[\overbar{\nu}](\overbar{\xi},\overbar{\eta})
		= 2\overbar{g}_{\overbar{\theta}}(\overbar{Z},\overbar{\nu}),
	\end{equation*}
	where
	\begin{multline*}
		\overbar{Z}=
		((\MAT{I}_p+\MAT{U}\MAT{U}^H)\herm(\MAT{\eta}_{\MAT{U}}\MAT{\xi}_{\MAT{U}}^H)\MAT{U} - \frac12\MAT{U}\herm(\MAT{\eta}_{\MAT{U}}^H\MAT{\xi}_{\MAT{U}}),
		\\
		-\MAT{\xi}_{\MAT{\Sigma}}\MAT{\Sigma}^{-1}\MAT{\eta}_{\MAT{\Sigma}}).
	\end{multline*}
	Since $\overbar{\nu}\in T_{\overbar{\theta}}\overbar{\mathcal{M}}_{p,k}$ and the projection map~\eqref{eq:rg_hlr:proj_product} is orthogonal according to~\eqref{eq:rg_hlr:metric}, projecting $\overbar{Z}$ on $T_{\overbar{\theta}}\overbar{\mathcal{M}}_{p,k}$ does not change the metric, \emph{i.e.}, $\overbar{g}_{\overbar{\theta}}(\overbar{Z},\overbar{\nu})=\overbar{g}_{\overbar{\theta}}(P_{\overbar{\theta}}(\overbar{Z}),\overbar{\nu})$.
	Thus,
	\begin{multline*}
		\D\overbar{g}_{\overbar{\theta}}[\overbar{\xi}](\overbar{\eta},\overbar{\nu})
		+ \D\overbar{g}_{\overbar{\theta}}[\overbar{\eta}](\overbar{\xi},\overbar{\nu})
		- \D\overbar{g}_{\overbar{\theta}}[\overbar{\nu}](\overbar{\xi},\overbar{\eta})
		= \\
		2\overbar{g}_{\overbar{\theta}}((\MAT{I}_p-\MAT{U}\MAT{U}^H)\herm(\MAT{\eta}_{\MAT{U}}\MAT{\xi}_{\MAT{U}}^H)\MAT{U}, -\herm(\MAT{\eta}_{\MAT{\Sigma}}\MAT{\Sigma}^{-1}\MAT{\xi}_{\MAT{\Sigma}}),\overbar{\nu}).
	\end{multline*}
	The same way, $\overbar{g}_{\overbar{\theta}}(\D\overbar{\eta}[\overbar{\xi}],\overbar{\nu})=\overbar{g}_{\overbar{\theta}}(P_{\overbar{\theta}}(\D\overbar{\eta}[\overbar{\xi}]),\overbar{\nu})$.
	Injecting these results in the Koszul formula, the Levi-Civita connection $\overbar{\nabla}_{\overbar{\xi}}\,\overbar{\eta}$ on $\overbar{\mathcal{M}}_{p,k}$ is finally obtained by identification.
	The Levi-Civita connection $\nabla_{\xi}\, \eta$ on $\mathcal{M}_{p,k}$ is then simply given by~\cite[proposition 5.3.3]{AMS08}.
\end{IEEEproof}

\vspace*{.5cm}
The geodesics in $\mathcal{M}_{p,k}$ associated with the metric induced by~\eqref{eq:rg_hlr:metric}, which generalize the concept of straight lines in a manifold, are given in proposition~\ref{prop:rg_hlr:geodesics}.
These geodesics are used to define a retraction on $\mathcal{M}_{p,k}$, \emph{i.e.}, a map from the tangent spaces back onto the manifold.
Unfortunately, an analytical formula for the geodesic between two points $\theta$ and $\widehat{\theta}$ in $\mathcal{M}_{p,k}$ is not known.
As a direct consequence, the Riemannian logarithm map and the Riemannian distance function on $\mathcal{M}_{p,k}$ are not known in closed~form.
\begin{proposition}
\label{prop:rg_hlr:geodesics}
	Let $\theta=\pi(\overbar{\theta})\in\mathcal{M}_{p,k}$ and $\xi=\D\pi(\overbar{\theta})[\overbar{\xi}]\in T_{\theta}\mathcal{M}_{p,k}$, where $\overbar{\xi}\in\mathcal{H}_{\overbar{\theta}}$.
	The representative in $\overbar{\mathcal{M}}_{p,k}$ of the geodesic in $\mathcal{M}_{p,k}$ associated with the metric induced by~\eqref{eq:rg_hlr:metric} starting at $\theta$ in the direction $\xi$ is%
	\footnote{
		The considered geodesic $\MAT{U}(t)$ on $\textup{St}_{p,k}$ is optimal (from a dimensionality point of view) only if $k\leq p/2$.
		If $k>{p}/2$, it is more advantageous to replace $\MAT{Q}$ with $\MAT{U}_{\perp}$ and $\MAT{R}$ with $\MAT{U}_{\perp}^H\MAT{\xi}_{\MAT{U}}$, where $\MAT{U}_{\perp}\in\textup{St}_{p,p-k}$ such that $\MAT{U}^H\MAT{U}_{\perp}=\MAT{0}$; see~\cite{EAS98}.
	}
	\begin{multline*}
		\overbar{\gamma}(t) = (\MAT{U}(t),\MAT{\Sigma}(t)) =
		\\
		\left([\MAT{U} \, \MAT{Q}] \exp t\begin{pmatrix}\MAT{U}^H\MAT{\xi}_{\MAT{U}} & -\MAT{R}^H \\ \MAT{R} & \MAT{0}\end{pmatrix} \begin{bmatrix}\MAT{I}_k\\\MAT{0}\end{bmatrix}, \right.
		\\
		\left. \MAT{\Sigma}^{\nicefrac12}\exp(t\MAT{\Sigma}^{\nicefrac{-1}2}\MAT{\xi}_{\MAT{\Sigma}}\MAT{\Sigma}^{\nicefrac{-1}2})\MAT{\Sigma}^{\nicefrac12} \right),
	\end{multline*}
	where $\MAT{Q}$ and $\MAT{R}$ correspond to the QR decomposition of $(\MAT{I}_p - \MAT{U}\MAT{U}^H)\MAT{\xi}_{\MAT{U}}$. 
\end{proposition}
\begin{IEEEproof}
	A direct proof that $\overbar{\gamma}(t)$ is a geodesic in $\overbar{\mathcal{M}}_{p,k}$ consists in verifying that it is solution of the differential equation $\overbar{\nabla}_{\dot{\overbar{\gamma}}(t)}\,\dot{\overbar{\gamma}}(t)=\MAT{0}$, where $\dot{\overbar{\gamma}}(t)$ is the derivative of $\overbar{\gamma}(t)$.
	However, it is enough to argue that $\MAT{U}(t)$ corresponds to the geodesic in $\textup{St}_{p,k}$ equipped with its canonical metric~\cite{EAS98} and $\MAT{\Sigma}(t)$ is the geodesic in $\mathcal{H}^{++}_k$ equipped with the considered affine invariant metric; see \emph{e.g.},~\cite{BGRB18}.

	To show that $\overbar{\gamma}(t)$ is a proper representative of the geodesic in $\mathcal{M}_{p,k}$, as $\pi$ is a Riemannian submersion, it suffices to show that $\overbar{\gamma}(t)$ stays horizontal in $\overbar{\mathcal{M}}_{p,k}$, \emph{i.e.}, $\dot{\overbar{\gamma}}(t)\in\mathcal{H}_{\overbar{\gamma}(t)}$~\cite[proposition 2.109]{GHL04}.
	One can check that $\MAT{U}(t)^H\dot{\MAT{U}}(t)=\MAT{U}^H\MAT{\xi}_{\MAT{U}}$, $\MAT{\Sigma}(t)^{-1}\dot{\MAT{\Sigma}}(t)=\MAT{\Sigma}^{-1}\MAT{\xi}_{\MAT{\Sigma}}$ and $\dot{\MAT{\Sigma}}(t)\MAT{\Sigma}(t)^{-1}=\MAT{\xi}_{\MAT{\Sigma}}\MAT{\Sigma}^{-1}$, which is enough to conclude. 
\end{IEEEproof}
\vspace*{.5cm}


\section{Riemannian optimization for robust covariance estimation}
\label{sec:ro_hlr}

We build a Riemannian optimization framework on $\mathcal{M}_{p,k}$ for robust estimation of covariance matrices admitting the structure~\eqref{eq:model:spiked_model}.
In section~\ref{subsec:ro_hlr:ro}, we provide the objects required to perform Riemannian optimization~\cite{AMS08} on $\mathcal{M}_{p,k}$, \emph{i.e.}, the Riemannian gradient and Hessian and a retraction, which corresponds to a second-order approximation of the geodesics of proposition~\ref{prop:rg_hlr:geodesics}.
In section~\ref{subsec:ro_hlr:obj_fun}, we develop tools to treat the family of cost functions of interest, which are originally defined on $\mathcal{H}^{++}_p$.
In particular, we deal with Tyler's $M$-estimator cost function defined in~\eqref{eq:model:tyler_HLR}.

\subsection{Riemannian optimization on $\mathcal{M}_{p,k}$}
\label{subsec:ro_hlr:ro}

Let $\overbar{f}:\overbar{\mathcal{M}}_{p,k}\to\mathbb{R}$ be an objective function that induces a function $f$ on the quotient $\mathcal{M}_{p,k}$, \emph{i.e.}, $\overbar{f}$ is invariant along the equivalence classes~\eqref{eq:model:eq_class}: for all $\overbar{\theta}\in\overbar{\mathcal{M}}_{p,k}$ and $\MAT{O}\in\mathcal{U}_k$, $\overbar{f}(\overbar{\theta})=\overbar{f}(\phi_{\MAT{O}}(\overbar{\theta}))$, where $\phi_{\MAT{O}}(\overbar{\theta})=(\MAT{U}\MAT{O},\MAT{O}^H\MAT{\Sigma}\MAT{O})$, as in section~\ref{sec:rg_hlr}.
To perform Riemannian optimization, it remains to define the Riemannian gradient and Hessian of $f$ along with a retraction on $\mathcal{M}_{p,k}$.
Proposition~\ref{prop:rgo_lr:grad_hess} provides formulas to compute the Riemannian gradient and Hessian of $f$ on $\mathcal{M}_{p,k}$ from the Euclidean gradient and Hessian of $\overbar{f}$ on $\overbar{\mathcal{M}}_{p,k}$.

\begin{proposition}
\label{prop:rgo_lr:grad_hess}
	Given $\theta=\pi(\overbar{\theta})\in\mathcal{M}_{p,k}$, the representative in $\mathcal{H}_{\overbar{\theta}}$ of the Riemannian gradient of $f$ at $\theta$ is the Riemannian gradient of $\overbar{f}$ at $\overbar{\theta}$, which is
	\begin{multline*}
		\grad_{\overbar{\mathcal{M}}_{p,k}}\overbar{f}(\overbar{\theta}) = 
		\left( \MAT{G}_{\MAT{U}}-\MAT{U}\MAT{G}_{\MAT{U}}^H\MAT{U}, \right.
		\\
		\left. \frac{\MAT{\Sigma}\herm(\MAT{G}_{\MAT{\Sigma}})\MAT{\Sigma}}{\alpha} - \frac{\beta\tr(\MAT{G}_{\MAT{\Sigma}}\MAT{\Sigma})}{\alpha(\alpha+k\beta)}\MAT{\Sigma} \right),
	\end{multline*}
	where $\grad_{\mathcal{E}}\overbar{f}(\overbar{\theta})=(\MAT{G}_{\MAT{U}},\MAT{G}_{\MAT{\Sigma}})$ is the Euclidean gradient of $\overbar{f}$ in $\mathbb{C}^{p\times k}\times\mathbb{C}^{k\times k}$.

	Given $\xi=\D\pi(\overbar{\theta})[\overbar{\xi}]\in T_{\theta}\mathcal{M}_{p,k}$, the representative in $\mathcal{H}_{\overbar{\theta}}$ of the Riemannian Hessian $\Hess_{\mathcal{M}_{p,k}}f(\theta)[\xi]$ of $f$ at $\theta$ in direction $\xi$ is
	\begin{equation*}
		\overbar{\Hess_{\mathcal{M}_{p,k}}f(\theta)[\xi]} = P^{\mathcal{H}}_{\overbar{\theta}}(\Hess_{\overbar{\mathcal{M}}_{p,k}}\overbar{f}(\overbar{\theta})[\overbar{\xi}]),
	\end{equation*}
	where $\Hess_{\overbar{\mathcal{M}}_{p,k}}\overbar{f}(\overbar{\theta})[\overbar{\xi}]$ is the Riemannian Hessian of $\overbar{f}$ at $\overbar{\theta}$ in direction $\overbar{\xi}$, given by
	\begin{multline*}
		\Hess_{\overbar{\mathcal{M}}_{p,k}}\overbar{f}(\overbar{\theta})[\overbar{\xi}] =
		\left( \MAT{H}_{\MAT{U}} - \MAT{U}\MAT{H}_{\MAT{U}}^H\MAT{U} - \MAT{U}\skewh(\MAT{G}_{\MAT{U}}^H\MAT{\xi}_{\MAT{U}})
		\right.
		\\
		- \skewh(\MAT{G}_{\MAT{U}}\MAT{\xi}_{\MAT{U}}^H)\MAT{U}
		- \frac12(\MAT{I}_p - \MAT{U}\MAT{U}^H)\MAT{\xi}_{\MAT{U}}\MAT{U}^H\MAT{G}_{\MAT{U}},
		\\
		\frac1{\alpha}(\MAT{\Sigma}\herm(\MAT{H}_{\MAT{\Sigma}})\MAT{\Sigma} + \herm(\MAT{\Sigma}\herm(\MAT{G}_{\MAT{\Sigma}})\MAT{\xi}_{\MAT{\Sigma}}))
		\\
		\left. - \frac{\beta\tr(\MAT{H}_{\MAT{\Sigma}}\MAT{\Sigma}+\MAT{G}_{\MAT{\Sigma}}\MAT{\xi}_{\MAT{\Sigma}})}{\alpha(\alpha+k\beta)}\MAT{\Sigma}\right),
	\end{multline*}
	where $\skewh$ returns the skew-Hermitian part of its argument and $\Hess_{\mathcal{E}}\overbar{f}(\overbar{\theta})[\overbar{\xi}]=(\MAT{H}_{\MAT{U}},\MAT{H}_{\MAT{\Sigma}})$ is the Euclidean Hessian of $\overbar{f}$ at $\overbar{\theta}$ in direction $\overbar{\xi}$, {i.e.,} $\Hess_{\mathcal{E}}\overbar{f}(\overbar{\theta})[\overbar{\xi}]=\D\grad_{\mathcal{E}}\overbar{f}(\overbar{\theta})[\overbar{\xi}]$.
\end{proposition}
\begin{IEEEproof}
	The Riemannian and Euclidean gradients of $\overbar{f}$ at $\overbar{\theta}$ are defined by
	\begin{equation*}
		\D\overbar{f}(\overbar{\theta})[\overbar{\xi}]
		= \langle \grad_{\overbar{\mathcal{M}}_{p,k}}\overbar{f}(\overbar{\theta}), \overbar{\xi} \rangle_{\overbar{\theta}} 
		= \langle \grad_{\mathcal{E}}\overbar{f}(\overbar{\theta}), \overbar{\xi} \rangle^{\mathcal{E}},
	\end{equation*}
	where $\langle\cdot,\cdot\rangle^{\mathcal{E}}$ is the Euclidean metric on $\mathbb{C}^{p\times k}\times\mathbb{C}^{k\times k}$, which is given by
	\begin{equation}
		\langle \xi, \eta \rangle^{\mathcal{E}} = \reel(\tr(\MAT{\xi}_{\MAT{U}}^H\MAT{\eta}_{\MAT{U}}) + \tr(\MAT{\xi}_{\MAT{\Sigma}}^H\MAT{\eta}_{\MAT{\Sigma}})).
	\label{eq:ro_hlr:metric_eucl_prod}
	\end{equation}
	Injecting the proposed formula for the gradient $\grad_{\overbar{\mathcal{M}}_{p,k}}\overbar{f}(\overbar{\theta})$ in the metric~\eqref{eq:rg_hlr:metric} shows that $\langle \grad_{\overbar{\mathcal{M}}_{p,k}}\overbar{f}(\overbar{\theta}), \overbar{\xi} \rangle_{\overbar{\theta}} $ is equal to $\langle \grad_{\mathcal{E}}\overbar{f}(\overbar{\theta}), \overbar{\xi} \rangle^{\mathcal{E}}$.
	To show that it is the Riemannian gradient of $\overbar{f}$ at $\overbar{\theta}\in\overbar{\mathcal{M}}_{p,k}$, we also need to check that it belongs to $T_{\overbar{\theta}}\overbar{\mathcal{M}}_{p,k}$ defined in~\eqref{eq:rg_hlr:tangent_product}, which is achieved with basic calculations.
	From~\cite{AMS08}, we further know that it belongs to $\mathcal{H}_{\overbar{\theta}}$ and that it is the representative of the Riemannian gradient of $f$ at $\theta\in\mathcal{M}_{p,k}$.

	The Riemannian Hessian of $\overbar{f}$ at $\overbar{\theta}$ in direction $\overbar{\xi}$ is defined as $\Hess_{\overbar{\mathcal{M}}_{p,k}}\overbar{f}(\overbar{\theta})[\overbar{\xi}]=\overbar{\nabla}_{\overbar{\xi}}\,\grad_{\overbar{\mathcal{M}}_{p,k}}\overbar{f}(\overbar{\theta})$~\cite{AMS08}.
	The result is obtained by plugging the formula of the gradient in the one of the Levi-Civita connection $\overbar{\nabla}$ on $\overbar{\mathcal{M}}_{p,k}$ defined in proposition~\ref{prop:rg_hlr:LC}.
	Finally, the representative of the Riemannian Hessian of $f$ at $\theta=\pi(\overbar{\theta})$ in direction $\xi=\D\pi(\overbar{\theta})[\overbar{\xi}]$ is obtained by definition of the Levi-Civita connection $\nabla$ on~$\mathcal{M}_{p,k}$, given in proposition~\ref{prop:rg_hlr:LC}.
\end{IEEEproof}

From the Riemannian gradient and Hessian, one can obtain a representative of a descent direction of $f$ at $\theta=\pi(\overbar{\theta})\in\mathcal{M}_{p,k}$ by selecting $\xi=\D\pi(\overbar{\theta})[\overbar{\xi}]\in T_{\theta}\mathcal{M}_{p,k}$ satisfying $\langle \grad_{\overbar{\mathcal{M}}_{p,k}}\overbar{f}(\overbar{\theta}), \overbar{\xi} \rangle_{\overbar{\theta}} < 0$.
A new point on the manifold is then achieved by a retraction on $\mathcal{M}_{p,k}$.
A natural choice is to take the Riemannian exponential map defined through the geodesics of proposition~\ref{prop:rg_hlr:geodesics}.
However, for numerical stability reasons, we rather choose a second order approximation of this exponential map, which, for $\theta=\pi(\overbar{\theta})\in\mathcal{M}_{p,k}$ and $\xi=\D\pi(\overbar{\theta})[\overbar{\xi}]\in T_{\theta}\mathcal{M}_{p,k}$, is represented by
\begin{multline}
	\overbar{R}_{\overbar{\theta}}(\overbar{\xi}) = \left( [\MAT{U} \, \MAT{Q}] \uf\circ\Gamma \begin{pmatrix}\MAT{U}^H\MAT{\xi}_{\MAT{U}} & -\MAT{R}^H \\ \MAT{R} & \MAT{0}\end{pmatrix} \begin{bmatrix}\MAT{I}_k\\\MAT{0}\end{bmatrix} \right.
	\\
	\left. \MAT{\Sigma}^{\nicefrac12}\Gamma(\MAT{\Sigma}^{\nicefrac{-1}2}\MAT{\xi}_{\MAT{\Sigma}}\MAT{\Sigma}^{\nicefrac{-1}2})\MAT{\Sigma}^{\nicefrac12} \right),
\label{eq:rgo_lr:retr}
\end{multline}
where $\uf$ returns the orthogonal factor of the polar decomposition and $\Gamma(\MAT{X})=\MAT{I}+\MAT{X}+\frac12\MAT{X}^2$ is a second order approximation of the matrix exponential.


With the tools developed in this section (and in section~\ref{sec:rg_hlr}), a large panel of first and second order Riemannian optimization algorithms can be employed to solve optimization problems on $\mathcal{M}_{p,k}$, such as gradient descent, conjugate gradient, BFGS, trust region, Newton, \emph{etc.}; see~\cite{AMS08} for details.
For example, given iterate $\theta_i=\pi(\overbar{\theta}_i)$, the Riemannian gradient descent algorithm yields iterate $\theta_{i+1}=\pi(\overbar{\theta}_{i+1})$ as
\begin{equation*}
	\theta_{i+1}=\pi\left(\overbar{R}_{\overbar{\theta}_i}\left( -t_i\grad_{\overbar{\mathcal{M}}_{p,k}}\overbar{f}(\overbar{\theta}_i) \right)\right),
\end{equation*}
where $t_i$ is the stepsize, which can for instance be computed with a line search~\cite{AMS08}.

\subsection{Robust covariance estimation}
\label{subsec:ro_hlr:obj_fun}

As detailed in section~\ref{sec:model}, we aim at estimating covariance matrices admitting the structure $\MAT{R}=\MAT{I}_p+\varphi(\theta)$, where $\varphi(\theta)=\overbar{\varphi}(\overbar{\theta})$, which is defined in~\eqref{eq:model:mapping_product2HLR}.
To that end, we are interested in objective functions $L:\mathcal{M}_{p,k}\to\mathbb{R}$ which have the form
\begin{equation}
	L(\theta) = L^{++}(\MAT{I}_p+\varphi(\theta)),
\label{eq:ro_hlr:cost_HPD2hlr}
\end{equation}
where $L^{++}:\mathcal{H}^{++}_p\to\mathbb{R}$ corresponds to an objective function for robust covariance estimation on $\mathcal{H}^{++}_p$, such as Tyler's $M$-estimator cost function~\eqref{eq:model:tyler_HPD}.
To perform Riemannian optimization of $L$ with the tools developed in section~\ref{subsec:ro_hlr:ro}, we simply need to have the Euclidean gradient and Hessian of $\overbar{L}=L\circ\pi$.
Proposition~\ref{prop:ro_hlr:grad_hess_HPD2hlr} shows that they can be obtained from those of~$L^{++}$.
For the Hessian, we need the directional derivative of $\overbar{\varphi}$ at $\overbar{\theta}$, which is given, for all $\overbar{\xi}\in T_{\overbar{\theta}}\overbar{\mathcal{M}}_{p,k}$, by
\begin{equation}
	\D\overbar{\varphi}(\overbar{\theta})[\overbar{\xi}]=\MAT{U}\MAT{\Sigma}\MAT{\xi}_{\MAT{U}}^H + \MAT{\xi}_{\MAT{U}}\MAT{\Sigma}\MAT{U}^H + \MAT{U}\MAT{\xi}_{\MAT{\Sigma}}\MAT{U}^H.
\label{eq:ro_hlr:Dmapping_product2HLR}
\end{equation}

\begin{proposition}
\label{prop:ro_hlr:grad_hess_HPD2hlr}
	The Euclidean gradient of $\overbar{L}=L\circ\pi$ at $\overbar{\theta}\in\overbar{\mathcal{M}}_{p,k}$ is given by
	\begin{equation*}
		\grad_{\mathcal{E}}\overbar{L}(\overbar{\theta}) = 
		(2\MAT{G}^{++}_{\overbar{\theta}}\MAT{U}\MAT{\Sigma},
		\MAT{U}^H\MAT{G}^{++}_{\overbar{\theta}}\MAT{U}),
	\end{equation*}
	where $\MAT{G}^{++}_{\overbar{\theta}}=\grad_{\mathfrak{E}}L^{++}(\MAT{I}_p+\overbar{\varphi}(\overbar{\theta}))$ is the Euclidean gradient of $L^{++}$ at $\MAT{I}_p+\overbar{\varphi}(\overbar{\theta})\in\mathcal{H}^{++}_p$, with $\overbar{\varphi}(\overbar{\theta})$ defined in~\eqref{eq:model:mapping_product2HLR}.

	The Euclidean Hessian of $\overbar{L}$ at $\overbar{\theta}$ in direction $\overbar{\xi}$ is
	\begin{multline*}
		\Hess_{\mathcal{E}}\overbar{L}(\overbar{\theta})[\overbar{\xi}]=
		(2\MAT{H}^{++}_{\overbar{\theta}}\MAT{U}\MAT{\Sigma} 
		+ 2\MAT{G}^{++}_{\overbar{\theta}}(\MAT{\xi}_{\MAT{U}}\MAT{\Sigma}+\MAT{U}\MAT{\xi}_{\MAT{\Sigma}}),
		\\
		\MAT{U}^H\MAT{H}^{++}_{\overbar{\theta}}\MAT{U}
		+ \MAT{U}^H\MAT{G}^{++}_{\overbar{\theta}}\MAT{\xi}_{\MAT{U}} + \MAT{\xi}_{\MAT{U}}^H\MAT{G}^{++}_{\overbar{\theta}}\MAT{U}
		),
	\end{multline*}
	where $\MAT{H}^{++}_{\overbar{\theta}}=\Hess_{\mathfrak{E}}L^{++}(\MAT{I}_p+\overbar{\varphi}(\overbar{\theta}))[\D\overbar{\varphi}(\overbar{\theta})[\overbar{\xi}]]$ is the Euclidean Hessian of $L^{++}$ at $\MAT{I}_p+\overbar{\varphi}(\overbar{\theta})\in\mathcal{H}^{++}_p$ in direction $\D\overbar{\varphi}(\overbar{\theta})[\overbar{\xi}]\in\mathcal{H}_p$, which is defined in~\eqref{eq:ro_hlr:Dmapping_product2HLR}.
\end{proposition}
\begin{IEEEproof}
	Let $\grad_{\mathcal{E}}\overbar{L}(\overbar{\theta})=(\MAT{G}_{\MAT{U}},\MAT{G}_{\MAT{\Sigma}})$.
	By definition,
	\begin{equation*}
		\D\overbar{f}(\overbar{\theta})[\overbar{\xi}] =
		\langle \grad_{\mathcal{E}}\overbar{L}(\overbar{\theta}), \overbar{\xi} \rangle^{\mathcal{E}},
	\end{equation*}
	where $\langle\cdot,\cdot\rangle^{\mathcal{E}}$ is defined in~\eqref{eq:ro_hlr:metric_eucl_prod}.
	We also have
	\begin{equation*}
		\begin{array}{rcl}
			\D\overbar{f}(\overbar{\theta})[\overbar{\xi}] & = & \D f^{++}(\MAT{I}_p + \overbar{\varphi}(\overbar{\theta}))[\D\overbar{\varphi}(\overbar{\theta})[\overbar{\xi}]]
			\\
			& = & \langle \MAT{G}^{++}_{\overbar{\theta}}, \D\overbar{\varphi}(\overbar{\theta})[\overbar{\xi}] \rangle^{\mathfrak{E}},
		\end{array}
	\end{equation*}
	where $\langle\cdot,\cdot\rangle^{\mathfrak{E}}$ is the Euclidean metric on $\mathbb{C}^{p\times p}$, which is
	\begin{equation*}
		\langle \MAT{\xi},\MAT{\eta} \rangle^{\mathfrak{E}} = \reel(\tr(\MAT{\xi}^H\MAT{\eta})).
	\end{equation*}
	We thus need to show that
	\begin{equation*}
		\reel(\tr(\MAT{G}_{\MAT{U}}^H\MAT{\xi}_{\MAT{U}}) + \tr(\MAT{G}_{\MAT{\Sigma}}^H\MAT{\xi}_{\MAT{\Sigma}}))
		= \reel(\tr(\MAT{G}^{++ \, H}_{\overbar{\theta}}\D\overbar{\varphi}(\overbar{\theta})[\overbar{\xi}])).
	\end{equation*}
	It is achieved by plugging the proposed formula for the Euclidean gradient $\grad_{\mathcal{E}}\overbar{L}(\overbar{\theta})=(\MAT{G}_{\MAT{U}},\MAT{G}_{\MAT{\Sigma}})$ and the definition of $\D\overbar{\varphi}(\overbar{\theta})[\overbar{\xi}]$ provided in~\eqref{eq:ro_hlr:Dmapping_product2HLR}.
	The Hessian is defined as $\Hess_{\mathcal{E}}\overbar{L}(\overbar{\theta})[\overbar{\xi}]=\D\grad_{\mathcal{E}}\overbar{L}(\overbar{\theta})[\overbar{\xi}]$.
	The proposed formula follows from basic calculations.
\end{IEEEproof}

To be able to compute Tyler's $M$-estimator on $\mathcal{M}_{p,k}$ from minimizing $L_{\textup{T}}$ defined in~\eqref{eq:model:tyler_HLR}, it remains to give the Euclidean gradient and Hessian of $L^{++}_{\textup{T}}$ defined in~\eqref{eq:model:tyler_HPD}.
To do so, we define $\Psi:\mathcal{H}^{++}_p\to\mathcal{H}_p$ and its directional derivative as
\begin{equation*}
	\begin{array}{rcl}
		\Psi(\MAT{R}) & = & \sum_i \frac{\VEC{x}_i\VEC{x}_i^H}{\VEC{x}_i^H\MAT{R}^{-1}\VEC{x}_i},
		\\[5pt]
		\D\Psi(\MAT{R})[\MAT{\xi}_{\MAT{R}}] & = & \sum_i \frac{\VEC{x}_i^H\MAT{R}^{-1}\MAT{\xi}_{\MAT{R}}\MAT{R}^{-1}\VEC{x}_i}{(\VEC{x}_i^H\MAT{R}^{-1}\VEC{x}_i)^2}\VEC{x}_i\VEC{x}_i^H.
	\end{array}
\end{equation*}
It follows that the Euclidean gradient of $L^{++}_{\textup{T}}$ at $\MAT{R}\in\mathcal{H}^{++}_p$ is
\begin{equation}
	\grad_{\mathfrak{E}} L^{++}_{\textup{T}}(\MAT{R}) = \MAT{R}^{-1}( n\MAT{R} - p\Psi(\MAT{R}) )\MAT{R}^{-1},
\end{equation}
and the Euclidean Hessian of $L^{++}_{\textup{T}}$ at $\MAT{R}\in\mathcal{H}^{++}_p$ in direction $\MAT{\xi}_{\MAT{R}}\in\mathcal{H}_p$ is
\begin{multline}
	\Hess_{\mathfrak{E}} L^{++}_{\textup{T}}(\MAT{R})[\MAT{\xi}_{\MAT{R}}] = 
	2p\MAT{R}^{-1}\herm(\MAT{\xi}_{\MAT{R}}\MAT{R}^{-1}\Psi(\MAT{R}))\MAT{R}^{-1}
	\\
	-\MAT{R}^{-1}(p\D\Psi(\MAT{R})[\MAT{\xi}_{\MAT{R}}]+n\MAT{\xi}_{\MAT{R}})\MAT{R}^{-1}.
\end{multline}

\section{Intrinsic Cram\'er-Rao lower bounds for sipked complex elliptically symmetric distributions}
\label{sec:icrb}




The manifold $\mathcal{M}_{p,k}$ admits a geometrical structure, described in section~\ref{sec:rg_hlr}, which can be exploited to measure the error of an unbiased estimator $\widehat{\theta}$ of the true parameter $\theta$ in $\mathcal{M}_{p,k}$.
A so-called lower intrinsic Cram\'er-Rao bound of such geometrical error measure can be obtained by exploiting the framework of~\cite{S05,B13}.
In section~\ref{subsec:icrb:err}, we define two different error measures: the first one is built from a divergence on $\mathcal{M}_{p,k}$ inspired by the one proposed in~\cite{BS09}; and the second one, which measures the subspace estimation error, is constructed from the Riemannian distance on the Grassmann manifold $\mathcal{G}_{p,k}$.
In section~\ref{subsec:icrb:fim}, we compute the Fisher information matrix on $\mathcal{M}_{p,k}$ associated with the distribution with pdf~\eqref{eq:model:ces_pdf_HLR}, which is needed to obtain the intrinsic Cram\'er-Rao bound.
We also study its structure in order to be able to bound the two error measures that we consider.
Finally, in section~\ref{subsec:icrb:ineq}, the intrinsic Cram\'er-bound inequalities are given.

\subsection{Estimation error measure}
\label{subsec:icrb:err}

We define two different error measures for any unbiased estimator $\widehat{\theta}=\pi(\MAT{\widehat{U}},\MAT{\widehat{\Sigma}})$ of the true parameter $\theta=\pi(\MAT{U},\MAT{\Sigma})$.
The first one is obtained from a proposed divergence function on $\mathcal{M}_{p,k}$, which is inspired from the one of~\cite{BS09}.
The second one measures the error of subspace estimation.
It is obtained from the distance on the Grassmann manifold $\mathcal{G}_{p,k}$ between $\spann(\MAT{\widehat{U}})$ and $\spann(\MAT{U})$.

In the general framework of~\cite{S05,B13}, the error on $\mathcal{M}_{p,k}$ of the unbiased estimator $\widehat{\theta}$ of $\theta$ is measured with the squared Riemannian distance $\delta_{\mathcal{M}_{p,k}}^2$ on $\mathcal{M}_{p,k}$, \emph{i.e.}, $\err_{\theta}(\widehat{\theta})=\delta_{\mathcal{M}_{p,k}}^2(\theta,\widehat{\theta})$.
However, as explained in section~\ref{sec:rg_hlr}, the Riemannian distance $\delta_{\mathcal{M}_{p,k}}(\theta,\widehat{\theta})$ on $\mathcal{M}_{p,k}$, which is the length (according to the metric induced by~\eqref{eq:rg_hlr:metric}) of the geodesic $\gamma=\pi\circ\overbar{\gamma}$ connecting $\theta$ and $\widehat{\theta}$, is not analytically known.
To overcome this issue, we define a divergence $d_{\mathcal{M}_{p,k}}(\theta,\widehat{\theta})$ on $\mathcal{M}_{p,k}$, which corresponds to the length (according to the metric induced by~\eqref{eq:rg_hlr:metric}) of a non-minimal curve $\widetilde{\gamma}$ connecting $\theta$ and $\widehat{\theta}$.
The error is then defined as
\begin{equation}
	\err^{\mathcal{M}_{p,k}}_{\theta}(\widehat{\theta}) = d_{\mathcal{M}_{p,k}}(\theta,\widehat{\theta}).
\label{eq:icrb:err_hlr}
\end{equation}
Moreover, by construction of $d_{\mathcal{M}_{p,k}}(\theta,\widehat{\theta})$, we have
\begin{equation*}
	d_{\mathcal{M}_{p,k}}(\theta,\widehat{\theta}) \geq \delta_{\mathcal{M}_{p,k}}^2(\theta,\widehat{\theta}).
\end{equation*}
The chosen divergence $d_{\mathcal{M}_{p,k}}$, which is inspired by the one proposed in~\cite{BS09}, is given in proposition~\ref{prop:icrb:div}.

\begin{proposition}
\label{prop:icrb:div}
	The function $d_{\mathcal{M}_{p,k}}:\mathcal{M}_{p,k}\times\mathcal{M}_{p,k}\to\mathbb{R}^+$, defined, for $\theta=\pi(\MAT{U},\MAT{\Sigma})$ and $\widehat{\theta}=\pi(\MAT{\widehat{U}},\MAT{\widehat{\Sigma}})$, as
	\begin{multline*}
		d_{\mathcal{M}_{p,k}}(\theta,\widehat{\theta}) =
		\alpha\norm{\log(\MAT{\Sigma}^{\nicefrac{-1}2}\MAT{O}\MAT{\widehat{O}}^H\MAT{\widehat{\Sigma}}\MAT{\widehat{O}}\MAT{O}^H\MAT{\Sigma}^{\nicefrac{-1}2})}_F^2
		\\
		+ \beta(\log\det(\MAT{\Sigma}^{-1}\MAT{O}\MAT{\widehat{O}}^H\MAT{\widehat{\Sigma}}\MAT{\widehat{O}}\MAT{O}^H))^2
		+ \norm{\MAT{\Theta}}_F^2,
	\end{multline*}
	where $\MAT{O}$, $\MAT{\widehat{O}}$ and $\MAT{\Theta}$ correspond to the singular value decomposition $\MAT{U}^H\MAT{\widehat{U}}=\MAT{O}\cos(\MAT{\Theta})\MAT{\widehat{O}}^H$, is a symmetric divergence function on $\mathcal{M}_{p,k}$ obtained by measuring the squared length (according to the metric induced by~\eqref{eq:rg_hlr:metric}) of the curve $\widetilde{\gamma}(t)=\pi(\MAT{\widetilde{U}}(t),\MAT{\widetilde{\Sigma}}(t))$, such that $(\MAT{\widetilde{U}}(t),\MAT{\widetilde{\Sigma}}(t))$ is the geodesic on $\overbar{\mathcal{M}}_{p,k}$ defined as
	\begin{equation*}
		\begin{array}{rcl}
			\MAT{\widetilde{U}}(t) & = & \MAT{U}\MAT{O}\cos(t\MAT{\Theta})\MAT{O}^H 
			\\
			& & +(\MAT{I}_p-\MAT{U}\MAT{U}^H)\MAT{\widehat{U}}\MAT{\widehat{O}}(\sin(\MAT{\Theta}))^{\dagger}\sin(t\MAT{\Theta})\MAT{O}^H,
			\\
			\MAT{\widetilde{\Sigma}}(t) & = & \MAT{\Sigma}^{\nicefrac12} (\MAT{\Sigma}^{\nicefrac{-1}2}\MAT{O}\MAT{\widehat{O}}^H\MAT{\widehat{\Sigma}}\MAT{O}\MAT{\widehat{O}}^H\MAT{\Sigma}^{\nicefrac{-1}2})^t \MAT{\Sigma}^{\nicefrac12},
		\end{array}
	\end{equation*}
	where $\cdot^{\dagger}$ and $\cdot^t=\exp(t\log(\cdot))$ are Moore-Penrose pseudo-inverse and matrix power functions, respectively.
\end{proposition}
\begin{IEEEproof}
	We aim to show that $\widetilde{\gamma}(t)=\pi(\MAT{\widetilde{U}}(t),\MAT{\widetilde{\Sigma}}(t))$ is a curve in $\mathcal{M}_{p,k}$ constructed from a geodesic $(\MAT{\widetilde{U}}(t),\MAT{\widetilde{\Sigma}}(t))$ in $\overbar{\mathcal{M}}_{p,k}$ and that measuring its squared length according to the metric induced by~\eqref{eq:rg_hlr:metric} yields the divergence $d_{\mathcal{M}_{p,k}}$ on $\mathcal{M}_{p,k}$.
	The problem we encounter while working with the geodesics $\overbar{\gamma}$ of proposition~\ref{prop:rg_hlr:geodesics} emanating from the horizontal space $\mathcal{H}_{\overbar{\theta}}$ of proposition~\ref{prop:rg_hlr:horizontal} is that it is not known analytically which direction $\overbar{\xi}\in\mathcal{H}_{\overbar{\theta}}$ connects $\theta=\pi(\overbar{\theta})$ to $\widehat{\theta}$.
	To overcome this issue, \cite{BS09} proposes to consider the alternative horizontal space
	\begin{equation}
		\widetilde{\mathcal{H}}_{\overbar{\theta}} = \{ \overbar{\xi}\in T_{\overbar{\theta}}\overbar{\mathcal{M}}_{p,k}: \MAT{U}^H\MAT{\xi}_{\MAT{U}}=\MAT{0} \}.
	\label{eq:icrb:horizontal_alt}
	\end{equation}
	$\widetilde{\mathcal{H}}_{\overbar{\theta}}$ still provides proper representatives of the elements in $T_{\theta}\mathcal{M}_{p,k}$, \emph{i.e.}, given $\xi\in T_{\theta}\mathcal{M}_{p,k}$, there is a unique $\overbar{\xi}\in\widetilde{\mathcal{H}}_{\overbar{\theta}}$ such that $\xi=\D\pi(\overbar{\theta})[\overbar{\xi}]$.
	This horizontal space is advantageous because the geodesics in $\overbar{\mathcal{M}}_{p,k}$ emanating from it are well characterized: the part of the geodesics that concerns $\MAT{U}$ coincides with the geodesics of the Grassmann manifold $\mathcal{G}_{p,k}$ while the part that concerns $\MAT{\Sigma}$ does not change.
	However, as $\pi$ is no longer a Riemannian submersion in this case, the resulting curves $\widetilde{\gamma}$ are not geodesics in $\mathcal{M}_{p,k}$.
	Given $\MAT{U}$ and $\MAT{\widehat{U}}$, the geodesics in $\mathcal{G}_{p,k}$ connecting $\spann(\MAT{U})$ and $\spann(\MAT{\widehat{U}})$ is $\spann(\MAT{\widetilde{U}}(t))$, where $\MAT{\widetilde{U}}(t)$ is defined above; see \emph{e.g.}~\cite{BS09,AMS04}.
	Since $\MAT{\widetilde{U}}(0)=\MAT{U}$ and $\MAT{\widetilde{U}}(1)=\MAT{\widehat{U}}\MAT{\widehat{O}}\MAT{O}^H$, we shall have $\MAT{\widetilde{\Sigma}}(0)=\MAT{\Sigma}$ and $\MAT{\widetilde{\Sigma}}(1)=\MAT{O}\MAT{\widehat{O}}^H\MAT{\widehat{\Sigma}}\MAT{\widehat{O}}\MAT{O}^H$ in order to obtain $\theta=\pi(\MAT{\widetilde{U}}(0),\MAT{\widetilde{\Sigma}}(0))$ and $\widehat{\theta}=\pi(\MAT{\widetilde{U}}(1),\MAT{\widetilde{\Sigma}}(1))$.
	It follows that $\MAT{\widetilde{\Sigma}}(t)$ is the geodesic on $\mathcal{H}^{++}_k$ defined as above.
	
	Finally, the squared length of $(\MAT{\widetilde{U}}(t),\MAT{\widetilde{\Sigma}}(t))$ in $\overbar{\mathcal{M}}_{p,k}$ according to metric~\eqref{eq:rg_hlr:metric} is the sum of the squared length of $\MAT{\widetilde{U}}(t)$ in $\textup{St}_{p,k}$ equiped with its canonical metric and of the squared length of $\MAT{\widetilde{\Sigma}}(t)$ in $\mathcal{H}^{++}_k$ equiped with the affine invariant metric.
	By construction, the squared length of $\MAT{\widetilde{U}}(t)$ corresponds to the squared Riemannian distance on $\mathcal{G}_{p,k}$~\cite{AMS04} between $\spann(\MAT{U})$ and $\spann(\MAT{\widehat{U}})$ and the one of $\MAT{\widetilde{\Sigma}}(t)$ is the squared Riemannian distance on $\mathcal{H}^{++}_k$~\cite{BGRB18} between $\MAT{\Sigma}$ and $\MAT{O}\MAT{\widehat{O}}^H\MAT{\widehat{\Sigma}}\MAT{\widehat{O}}\MAT{O}^H$.
	We thus obtain the proposed function $d_{\mathcal{M}_{p,k}}(\theta,\widehat{\theta})$ on $\mathcal{M}_{p,k}$.
	Furthermore, it is readily checked that it is a symmetric divergence function on $\mathcal{M}_{p,k}$, \emph{i.e.}, $d_{\mathcal{M}_{p,k}}(\theta,\widehat{\theta})\geq0$, with equality if and only if $\theta=\widehat{\theta}$ and $d_{\mathcal{M}_{p,k}}(\theta,\widehat{\theta})=d_{\mathcal{M}_{p,k}}(\widehat{\theta},\theta)$.
\end{IEEEproof}

The principal subspace of $\varphi(\theta)\in\mathcal{H}^{+}_{p,k}$ is given by $\spann(\MAT{U})\in\mathcal{G}_{p,k}$, which is estimated by $\spann(\MAT{\widehat{U}})$.
To measure the error of subspace estimation, we exploit the Riemannian distance function on the Grassmann manifold $\mathcal{G}_{p,k}$ equiped with the Riemannian metric induced by the part of metric~\eqref{eq:rg_hlr:metric} that depends on the component $\MAT{U}$.
Thus, the subspace estimation error of the unbiased estimator $\widehat{\theta}=\pi(\MAT{\widehat{U}},\MAT{\widehat{\Sigma}})$ of the true parameter $\theta=\pi(\MAT{U},\MAT{\Sigma})$ is
\begin{equation}
	\err^{\mathcal{G}_{p,k}}_{\theta}(\widehat{\theta}) = \delta_{\mathcal{G}_{p,k}}^2(\spann(\MAT{U}),\spann(\MAT{\widehat{U}})),
\label{eq:icrb:err_Gr}
\end{equation}
where $\delta_{\mathcal{G}_{p,k}}$ is the Riemannian distance function on Grassmann (see \emph{e.g.}~\cite{AMS04}), which is given by
\begin{equation*}
	\delta_{\mathcal{G}_{p,k}}^2(\spann(\MAT{U}),\spann(\MAT{\widehat{U}})) = \norm{\MAT{\Theta}}_F^2,
\end{equation*}
with $\MAT{\Theta}$ defined as in proposition~\ref{prop:icrb:div}.

\subsection{Fisher information matrix}
\label{subsec:icrb:fim}

Given $\theta=\pi(\overbar{\theta})\in\mathcal{M}_{p,k}$, we compute the Fisher information matrix $\MAT{F}_{\overbar{\theta}}$ corresponding to the distribution with pdf~\eqref{eq:model:ces_pdf_HLR} and study its structure.
%
%
In order to obtain $\MAT{F}_{\overbar{\theta}}$, we need to define~\cite{B13}:
\emph{(i)} the Fisher information metric $g^{\overbar{\mathcal{M}}_{p,k}}_{\overbar{\theta}}$ on $\overbar{\mathcal{M}}_{p,k}$ of distribution with pdf~\eqref{eq:model:ces_pdf_HLR}
and
\emph{(ii)} an orthonormal basis $\{e_{\overbar{\theta}}^q\}_q$ of the tangent space $T_{\overbar{\theta}}\overbar{\mathcal{M}}_{p,k}$ according to metric~\eqref{eq:rg_hlr:metric}. 
The $q\ell^{\textup{th}}$ element of $\MAT{F}_{\overbar{\theta}}$ is then defined~as
\begin{equation}
	(\MAT{F}_{\overbar{\theta}})_{q\ell} = g^{\overbar{\mathcal{M}}_{p,k}}_{\overbar{\theta}}(e_{\overbar{\theta}}^q,e_{\overbar{\theta}}^\ell).
\label{eq:icrb:FIM}
\end{equation}
Notice that due to the invariance with respect to the action of unitary matrices in $\mathcal{U}_k$ described in section~\ref{subsec:model:cov_hlr}, $\MAT{F}_{\overbar{\theta}}$, whose size is $2pk\times 2pk$, has rank $2pk-k^2$.

Concerning the Fisher information metric, we first give the general proposition~\ref{prop:icrb:fim}, which allows to obtain the Fisher information on a manifold $\mathcal{M}$ from the one on another manifold $\mathcal{N}$ when the pdf on $\mathcal{M}$ is defined through the one on $\mathcal{N}$ and a smooth mapping $\psi:\mathcal{M}\to\mathcal{N}$.
The Fisher information metric on $\overbar{\mathcal{M}}_{p,k}$ corresponding to the pdf~\eqref{eq:model:ces_pdf_HLR} is then obtained from the one corresponding to the pdf~\eqref{eq:model:ces_pdf_HPD} on $\mathcal{H}^{++}_p$ in corollary~\ref{corol:icrb:fim}.

\begin{proposition}
\label{prop:icrb:fim}
	Let two manifolds $\mathcal{M}$, $\mathcal{N}$ and the smooth mapping $\psi:\mathcal{M}\to\mathcal{N}$.
	Consider the pdf on $\mathcal{M}$
	\begin{equation*}
		f^{\mathcal{M}}(x|\theta) = f^{\mathcal{N}}(x|\psi(\theta)),
	\end{equation*}
	where $f^{\mathcal{N}}$ is a pdf on $\mathcal{N}$ whose Fisher information metric is $g^{\mathcal{N}}$.
	It follows that the Fisher information metric $g^{\mathcal{M}}$ on $\mathcal{M}$ associated with $f^{\mathcal{M}}$ is, given $\theta\in\mathcal{M}$ and $\xi$, $\eta\in T_{\theta}\mathcal{M}$,
	\begin{equation*}
		g^{\mathcal{M}}_{\theta}(\xi,\eta) = g^{\mathcal{N}}_{\psi(\theta)}(\D\psi(\theta)[\xi],\D\psi(\theta)[\eta]).
	\end{equation*}
\end{proposition}
\begin{IEEEproof}
	Let $L^{\mathcal{M}}_x(\theta)$ the log-likelihood on $\mathcal{M}$ of $f^{\mathcal{M}}(x|\theta)$.
	By definition, $L^{\mathcal{M}}_x(\theta)=L^{\mathcal{N}}_x(\psi(\theta))$ and
	\begin{equation*}
		\begin{array}{l}
			g^{\mathcal{M}}_{\theta}(\xi,\eta)
			= \mathbb{E}_x\left[ \D L^{\mathcal{M}}_x(\theta)[\xi] \D L^{\mathcal{M}}_x(\theta)[\eta] \right]
			\\[5pt] \quad
			= \mathbb{E}_x\left[ \D L^{\mathcal{N}}_x(\psi(\theta))[\D\psi(\theta)[\xi]] \D L^{\mathcal{N}}_x(\psi(\theta))[\D\psi(\theta)[\eta]] \right]
			\\[5pt] \quad
			= g^{\mathcal{N}}_{\psi(\theta)}(\D\psi(\theta)[\xi],\D\psi(\theta)[\eta]).
		\end{array}
	\end{equation*}
\end{IEEEproof}

\begin{corollary}
\label{corol:icrb:fim}
	The Fisher information metric on $\overbar{\mathcal{M}}_{p,k}$ corresponding to the pdf~\eqref{eq:model:ces_pdf_HLR} is, for $\overbar{\theta}\in\overbar{\mathcal{M}}_{p,k}$ and $\overbar{\xi}$, $\overbar{\eta}\in T_{\overbar{\theta}}\overbar{\mathcal{M}}_{p,k}$,
	\begin{equation*}
		g^{\overbar{\mathcal{M}}_{p,k}}_{\overbar{\theta}}(\overbar{\xi},\overbar{\eta}) = 
		g^{\mathcal{H}^{++}_p}_{\MAT{I}_p + \overbar{\varphi}(\overbar{\theta})}(\D\overbar{\varphi}(\overbar{\theta})[\overbar{\xi}],\D\overbar{\varphi}(\overbar{\theta})[\overbar{\eta}]),
	\end{equation*}
	where $\overbar{\varphi}(\overbar{\theta})$ and $\D\overbar{\varphi}(\overbar{\theta})[\overbar{\xi}]$ are defined in~\eqref{eq:model:mapping_product2HLR} and~\eqref{eq:ro_hlr:Dmapping_product2HLR}, and
	\begin{multline*}
		g^{\mathcal{H}^{++}_p}_{\MAT{R}}(\MAT{\xi}_{\MAT{R}},\MAT{\eta}_{\MAT{R}}) =
		n\alpha^{++}\tr(\MAT{R}^{-1}\MAT{\xi}_{\MAT{R}}\MAT{R}^{-1}\MAT{\eta}_{\MAT{R}})
		\\
		+ n(\alpha^{++}-1)\tr(\MAT{R}^{-1}\MAT{\xi}_{\MAT{R}})\tr(\MAT{R}^{-1}\MAT{\xi}_{\MAT{R}})
	\end{multline*}
	is the Fisher information on $\mathcal{H}^{++}_p$ associated with the pdf~\eqref{eq:model:ces_pdf_HPD}, where $\alpha^{++}$ is a scalar that only depends on the density generator $g$ in~\eqref{eq:model:ces_pdf_HPD}~\cite{BGRB18}.
\end{corollary}

It remains to provide an orthonormal basis on the tangent space $T_{\overbar{\theta}}\overbar{\mathcal{M}}_{p,k}$ of $\overbar{\theta}\in\overbar{\mathcal{M}}_{p,k}$ according to metric~\eqref{eq:rg_hlr:metric} to be able to compute the Fisher information matrix $\MAT{F}_{\overbar{\theta}}$.
This is done in proposition~\ref{prop:icrb:basis}.

\begin{proposition}
\label{prop:icrb:basis}
	Given $\overbar{\theta}\in\overbar{\mathcal{M}}_{p,k}$, an orthonormal basis $\{e_{\overbar{\theta}}^q\}_{1\leq q\leq2pk}$ of the tangent space $T_{\overbar{\theta}}\overbar{\mathcal{M}}_{p,k}$ is given by
	\begin{multline*}
		\left\{
		\{(\MAT{e}_{\MAT{U}_{\perp}}^{ij}, \MAT{0}), (\MAT{\widetilde{e}}_{\MAT{U}_{\perp}}^{ij}, \MAT{0})\}_{\substack{1\leq i\leq p-k \\ 1\leq j\leq k}},
		\{(\MAT{e}_{\MAT{U}}^{ij}, \MAT{0})\}_{1\leq j<i\leq k},
		\right.
		\\
		\left.
		\{(\MAT{\widetilde{e}}_{\MAT{U}}^{ij}, \MAT{0})\}_{1\leq j\leq i\leq k},
		\{(\MAT{0}, \MAT{e}_{\MAT{\Sigma}}^{ij})\}_{1\leq j\leq i\leq k},
		\{(\MAT{0}, \MAT{\widetilde{e}}_{\MAT{\Sigma}}^{ij})\}_{1\leq j<i\leq k}
		\right\},
	\end{multline*}
	where 
	\begin{itemize}
		\item $\MAT{e}_{\MAT{U}_{\perp}}^{ij}=\MAT{U}_{\perp}\MAT{K}^{ij}$, $\MAT{\widetilde{e}}_{\MAT{U}_{\perp}}^{ij}=\I\MAT{U}_{\perp}\MAT{K}^{ij}$:
		$\MAT{U}_{\perp}\in\textup{St}_{p,p-k}$, $\MAT{U}^H\MAT{U}_{\perp}=\MAT{0}$;
		$\MAT{K}^{ij}\in\mathbb{R}^{(p-k)\times k}$, its $ij^{\textup{th}}$ element is $1$, zeros elsewhere.
		\item $\MAT{e}_{\MAT{U}}^{ij}=\MAT{U}\MAT{\Omega}^{ij}$:
		$\MAT{\Omega}^{ij}\in\mathcal{H}^{\perp}_k$, its $ij^{\textup{th}}$ and $ji^{\textup{th}}$ elements are $1$ and $-1$, zeros elsewhere.
		\item $\MAT{\widetilde{e}}_{\MAT{U}}^{ij}=\MAT{U}\MAT{\widetilde{\Omega}}^{ij}$:
		$\MAT{\widetilde{\Omega}}^{ii}\in\mathcal{H}^{\perp}_k$, its $ii^{\textup{th}}$ element is $\sqrt{2}\I$, zeros elsewhere;
		$\MAT{\widetilde{\Omega}}^{ij}\in\mathcal{H}^{\perp}_k$, $i>j$, its $ij^{\textup{th}}$ and $ji^{\textup{th}}$ elements are $\I$, zeros elsewhere.
		\item $\MAT{e}_{\MAT{\Sigma}}^{ij}=\frac1{\sqrt{\alpha}}\MAT{\Sigma}^{\nicefrac12}\MAT{H}^{ij}\MAT{\Sigma}^{\nicefrac12} + \frac{\sqrt{\alpha}-\sqrt{\alpha+k\beta}}{k\sqrt{\alpha}\sqrt{\alpha+k\beta}}\tr(\MAT{H}^{ij})\MAT{\Sigma}$:
		$\MAT{H}^{ii}\in\mathcal{H}_k$, its $ii^{\textup{th}}$ element is $1$, zeros elsewhere;
		$\MAT{H}^{ij}\in\mathcal{H}_k$, $i>j$, its $ij^{\textup{th}}$ and $ji^{\textup{th}}$ elements are $\nicefrac1{\sqrt{2}}$, zeros elsewhere.
		\item $\MAT{\widetilde{e}}_{\MAT{\Sigma}}^{ij}=\frac1{\sqrt{\alpha}}\MAT{\Sigma}^{\nicefrac12}\MAT{\widetilde{H}}^{ij}\MAT{\Sigma}^{\nicefrac12}$:
		$\MAT{\widetilde{H}}^{ij}\in\mathcal{H}_k$, its $ij^{\textup{th}}$ and $ji^{\textup{th}}$ elements are $\nicefrac{\I}{\sqrt{2}}$ and $\nicefrac{-\I}{\sqrt{2}}$, zeros elsewhere.
	\end{itemize}
\end{proposition}
\begin{IEEEproof}
	By definition, it suffices to check that for all $1\leq p,\ell\leq 2pk$, $p\neq\ell$, $\langle e_{\overbar{\theta}}^q, e_{\overbar{\theta}}^q\rangle_{\overbar{\theta}}=1$ and $\langle e_{\overbar{\theta}}^q, e_{\overbar{\theta}}^{\ell}\rangle_{\overbar{\theta}}=0$, which is achieved by basic calculations.
\end{IEEEproof}

In proposition~\ref{prop:icrb:FIM_struct}, we study the structure of the Fisher information matrix $\MAT{F}_{\overbar{\theta}}$ on $\overbar{\mathcal{M}}_{p,k}$ corresponding to the pdf~\eqref{eq:model:ces_pdf_HLR}.

\begin{proposition}
\label{prop:icrb:FIM_struct}
	The Fisher information matrix $\MAT{F}_{\overbar{\theta}}$ on $\overbar{\mathcal{M}}_{p,k}$ of the pdf~\eqref{eq:model:ces_pdf_HLR} admits the structure
	\begin{equation}
	\label{eq:structure_fim}
		\MAT{F}_{\overbar{\theta}} =
		\begin{pmatrix}
			\MAT{F}_{\MAT{U}_{\perp}} & \MAT{0}                        & \MAT{0}                        \\
			\MAT{0}                   & \MAT{F}_{\MAT{U}}              & \MAT{F}_{\MAT{U},\MAT{\Sigma}} \\
			\MAT{0}                   & \MAT{F}_{\MAT{\Sigma},\MAT{U}} & \MAT{F}_{\MAT{\Sigma}}
		\end{pmatrix},
	\end{equation}
	where $\MAT{F}_{\MAT{U}_{\perp}}\in\mathbb{R}^{2(p-k)k\times2(p-k)k}$ is the block obtained from the elements $\{(\MAT{e}_{\MAT{U}_{\perp}}^{ij}, \MAT{0}), (\MAT{\widetilde{e}}_{\MAT{U}_{\perp}}^{ij}, \MAT{0})\}$ of the orthonormal basis of $T_{\overbar{\theta}}\overbar{\mathcal{M}}_{p,k}$ given in proposition~\ref{prop:icrb:basis};
	and $\MAT{F}_{\MAT{U}}$, $\MAT{F}_{\MAT{U},\MAT{\Sigma}}$, $\MAT{F}_{\MAT{\Sigma},\MAT{U}}$, $\MAT{F}_{\MAT{\Sigma}}\in\mathbb{R}^{k^2\times k^2}$ are the blocks obtained from the remaining elements of the basis.
	Further notice that $\MAT{F}_{\MAT{U}_{\perp}}\in\mathbb{R}^{2(p-k)k\times2(p-k)k}$, $\MAT{F}_{\MAT{U}}\in\mathbb{R}^{k^2\times k^2}$ and $\MAT{F}_{\MAT{\Sigma}}\in\mathbb{R}^{k^2\times k^2}$ are of full rank, and
	\begin{equation*}
		\begin{pmatrix}
			\MAT{F}_{\MAT{U}}              & \MAT{F}_{\MAT{U},\MAT{\Sigma}} \\
			\MAT{F}_{\MAT{\Sigma},\MAT{U}} & \MAT{F}_{\MAT{\Sigma}}
		\end{pmatrix}
		\in \mathbb{R}^{2k^2\times2k^2}
	\end{equation*}
	has rank $k^2$.
\end{proposition}
\begin{IEEEproof}
	Every tangent vector $\MAT{\xi}_{\MAT{U}}\in T_{\MAT{U}}\textup{St}_{p,k}$ can be decomposed as $\MAT{\xi}_{\MAT{U}}=\MAT{U}\MAT{\Omega}_{\xi} + \MAT{U}_{\perp}\MAT{K}_{\xi}$, where $\MAT{U}_{\perp}\in\textup{St}_{p,p-k}$ such that $\MAT{U}^H\MAT{U}_{\perp}=\MAT{0}$, $\MAT{\Omega}_{\xi}\in\mathcal{H}^{\perp}_k$ and $\MAT{K}_{\xi}\in\mathbb{C}^{(p-k)\times k}$.
	Thus, $\overbar{\xi}\in T_{\overbar{\theta}}\overbar{\mathcal{M}}_{p,k}$ can be decomposed as
	\begin{equation*}
		\overbar{\xi}
		= \overbar{\xi}^{\MAT{U}} + \overbar{\xi}^{\MAT{U}_{\perp}} + \overbar{\xi}^{\MAT{\Sigma}}
		= (\MAT{U}\MAT{\Omega}_{\xi},\MAT{0}) + (\MAT{U}_{\perp}\MAT{K}_{\xi},\MAT{0}) + (\MAT{0},\MAT{\xi}_{\MAT{\Sigma}}).
	\end{equation*}
	By linearity of $g^{\overbar{\mathcal{M}}_{p,k}}_{\overbar{\theta}}$ defined in corollary~\ref{corol:icrb:fim}, we have
	\begin{equation*}
		\begin{array}{rcl}
		g^{\overbar{\mathcal{M}}_{p,k}}_{\overbar{\theta}}(\overbar{\xi},\overbar{\xi})
		& = & g^{\overbar{\mathcal{M}}_{p,k}}_{\overbar{\theta}}(\overbar{\xi}^{\MAT{U}_{\perp}},\overbar{\xi}^{\MAT{U}_{\perp}})
		+ g^{\overbar{\mathcal{M}}_{p,k}}_{\overbar{\theta}}(\overbar{\xi}^{\MAT{U}},\overbar{\xi}^{\MAT{U}})
		\\ & &
		+ g^{\overbar{\mathcal{M}}_{p,k}}_{\overbar{\theta}}(\overbar{\xi}^{\MAT{\Sigma}},\overbar{\xi}^{\MAT{\Sigma}})
		+ 2g^{\overbar{\mathcal{M}}_{p,k}}_{\overbar{\theta}}(\overbar{\xi}^{\MAT{U}},\overbar{\xi}^{\MAT{\Sigma}})
		\\ & &
		+ 2 g^{\overbar{\mathcal{M}}_{p,k}}_{\overbar{\theta}}(\overbar{\xi}^{\MAT{U}_{\perp}},\overbar{\xi}^{\MAT{U}})
		+ 2g^{\overbar{\mathcal{M}}_{p,k}}_{\overbar{\theta}}(\overbar{\xi}^{\MAT{U}_{\perp}},\overbar{\xi}^{\MAT{\Sigma}}).
		\end{array}
	\end{equation*}
	To show that $\MAT{F}_{\overbar{\theta}}$ has the proposed form, it suffices to prove that $g^{\overbar{\mathcal{M}}_{p,k}}_{\overbar{\theta}}(\overbar{\xi}^{\MAT{U}_{\perp}},\overbar{\xi}^{\MAT{U}})=g^{\overbar{\mathcal{M}}_{p,k}}_{\overbar{\theta}}(\overbar{\xi}^{\MAT{U}_{\perp}},\overbar{\xi}^{\MAT{\Sigma}})=0$.
	From~\eqref{eq:ro_hlr:Dmapping_product2HLR}, we obtain
	\begin{equation*}
		\begin{array}{rcl}
			\D\overbar{\varphi}(\overbar{\theta})[\overbar{\xi}^{\MAT{U}_{\perp}}] & = & \MAT{U}\MAT{\Sigma}\MAT{K}_{\xi}^H\MAT{U}_{\perp}^H + \MAT{U}_{\perp}\MAT{K}_{\xi}\MAT{\Sigma}\MAT{U}^H,
			\\[2pt]
			\D\overbar{\varphi}(\overbar{\theta})[\overbar{\xi}^{\MAT{U}}] & = & \MAT{U}\MAT{\Sigma}\MAT{\Omega}_{\xi}^H\MAT{U}^H + \MAT{U}\MAT{\Omega}_{\xi}\MAT{\Sigma}\MAT{U}^H,
			\\[1pt]
			\D\overbar{\varphi}(\overbar{\theta})[\overbar{\xi}^{\MAT{\Sigma}}] & = &  \MAT{U}\MAT{\xi}_{\MAT{\Sigma}}\MAT{U}^H.
		\end{array}
	\end{equation*}
	The Woodbury identity $\overbar{\varphi}(\overbar{\theta})^{-1}=\MAT{I}_p-\MAT{U}\MAT{\Xi}\MAT{U}^H$, where $\MAT{\Xi}=(\MAT{I}_k+\MAT{\Sigma }^{-1})^{-1}$, and $\MAT{U}^H\MAT{U}_{\perp}=\MAT{0}$ lead to
	\begin{multline*}
		\overbar{\varphi}(\overbar{\theta})^{-1}\D\overbar{\varphi}(\overbar{\theta})[\overbar{\xi}^{\MAT{U}_{\perp}}] =
		\MAT{U}\MAT{\Sigma}\MAT{K}_{\xi}^H\MAT{U}_{\perp}^H  
		+ \MAT{U}_{\perp}\MAT{K}_{\xi}\MAT{\Sigma}\MAT{U}^H
		\\
		- \MAT{U}\MAT{\Xi}\MAT{\Sigma}\MAT{K}_{\xi}^H\MAT{U}_{\perp}^H,
	\end{multline*}
	from which one can check that $\tr(\overbar{\varphi}(\overbar{\theta})^{-1}\D\overbar{\varphi}(\overbar{\theta})[\overbar{\xi}^{\MAT{U}_{\perp}}])=0$.
	Furthermore, the previous expression yields
	\begin{multline*}
		\overbar{\varphi}(\overbar{\theta})^{-1} \D\overbar{\varphi}(\overbar{\theta})[\overbar{\xi}^{\MAT{U}_{\perp}}] \overbar{\varphi}(\overbar{\theta})^{-1} =
		\MAT{U}_{\perp}\MAT{K}_{\xi}\MAT{\Sigma}(\MAT{I}_k - \MAT{\Xi})\MAT{U}^H
		\\
 		+ \MAT{U}(\MAT{I}_k - \MAT{\Xi})\MAT{\Sigma}\MAT{K}_{\xi}^H\MAT{U}_{\perp}^H.
	\end{multline*}
	From this, it is readily checked that $g^{\overbar{\mathcal{M}}_{p,k}}_{\overbar{\theta}}(\overbar{\xi}^{\MAT{U}_{\perp}},\overbar{\xi}^{\MAT{U}})=g^{\overbar{\mathcal{M}}_{p,k}}_{\overbar{\theta}}(\overbar{\xi}^{\MAT{U}_{\perp}},\overbar{\xi}^{\MAT{\Sigma}})=0$.

	Finally, to show that $\MAT{F}_{\MAT{U}_{\perp}}\in\mathbb{R}^{2(p-k)k\times2(p-k)k}$, $\MAT{F}_{\MAT{U}}\in\mathbb{R}^{k^2\times k^2}$ and $\MAT{F}_{\MAT{\Sigma}}\in\mathbb{R}^{k^2\times k^2}$ are of full rank, it is enough to verify that $\overbar{\xi}^{\MAT{U}_{\perp}}\mapsto g^{\overbar{\mathcal{M}}_{p,k}}_{\overbar{\theta}}(\overbar{\xi}^{\MAT{U}_{\perp}},\overbar{\xi}^{\MAT{U}_{\perp}})$, $\overbar{\xi}^{\MAT{U}}\mapsto g^{\overbar{\mathcal{M}}_{p,k}}_{\overbar{\theta}}(\overbar{\xi}^{\MAT{U}},\overbar{\xi}^{\MAT{U}})$ and $\overbar{\xi}^{\MAT{\Sigma}}\mapsto g^{\overbar{\mathcal{M}}_{p,k}}_{\overbar{\theta}}(\overbar{\xi}^{\MAT{\Sigma}},\overbar{\xi}^{\MAT{\Sigma}})$ are positive definite.
	The rank of
	\begin{equation*}
		\begin{pmatrix}
			\MAT{F}_{\MAT{U}}              & \MAT{F}_{\MAT{U},\MAT{\Sigma}} \\
			\MAT{F}_{\MAT{\Sigma},\MAT{U}} & \MAT{F}_{\MAT{\Sigma}}
		\end{pmatrix}
		\in \mathbb{R}^{2k^2\times2k^2}
	\end{equation*}
	is given by subtracting the rank of $\MAT{F}_{\MAT{U}_{\perp}}$ to the one of $\MAT{F}_{\overbar{\theta}}$.
\end{IEEEproof}

From~\cite{B13}, we know that the Fisher information matrix $\MAT{F}_{\overbar{\theta}}$ is adapted to the proposed geometry of the parameter manifold $\mathcal{M}_{p,k}$, \emph{i.e.}, it is well suited to bound the error measured through the distance function $\delta_{\mathcal{M}_{p,k}}$ on $\mathcal{M}_{p,k}$.
However, as this distance is not analytically known, the error is measured with the divergence of proposition~\ref{prop:icrb:div} in this work.
Recall that this divergence is obtained by considering the horizontal space $\widetilde{\mathcal{H}}_{\overbar{\theta}}$ defined in~\eqref{eq:icrb:horizontal_alt} instead of the one of proposition~\ref{prop:rg_hlr:horizontal}.
Hence, a Fisher information matrix $\MAT{\widetilde{F}}_{\overbar{\theta}}$ which appears well-suited to the divergence $d_{\mathcal{M}_{p,k}}$ is constructed by taking an orthonormal basis of $\widetilde{\mathcal{H}}_{\overbar{\theta}}$ according to metric~\eqref{eq:rg_hlr:metric}.
Such a basis is formed by the $2pk-k^2$ elements $\{(\MAT{e}_{\MAT{U}_{\perp}}^{ij}, \MAT{0}), (\MAT{\widetilde{e}}_{\MAT{U}_{\perp}}^{ij}, \MAT{0})\}$, $\{(\MAT{0}, \MAT{e}_{\MAT{\Sigma}}^{ij})\}$ and $\{(\MAT{0}, \MAT{\widetilde{e}}_{\MAT{\Sigma}}^{ij})\}$ defined in proposition~\ref{prop:icrb:basis}.
The resulting Fisher information matrix $\MAT{\widetilde{F}}_{\overbar{\theta}}$ is%
\footnote{
	Interestingly, $\MAT{\widetilde{F}}_{\overbar{\theta}}$ corresponds to the Fisher information matrix obtained in~\cite{S05} from a different reasonning for the Gaussian case ($\alpha^{++}=1$ in corollary~\ref{corol:icrb:fim}) and for $\alpha=1$ and $\beta=0$ in metric~\eqref{eq:rg_hlr:metric}.
}
\begin{equation}
	\MAT{\widetilde{F}}_{\overbar{\theta}} =
	\begin{pmatrix}
		\MAT{F}_{\MAT{U}_{\perp}} & \MAT{0}                \\
		\MAT{0}                   & \MAT{F}_{\MAT{\Sigma}}
	\end{pmatrix}.
\label{eq:icrb:FIM_alt}
\end{equation}
Its size is $(2pk-k^2)\times(2pk-k^2)$ and it has full rank.

\subsection{Inequalities}
\label{subsec:icrb:ineq}

Finally, we derive intrinsic Cram\'er-Rao lower bounds~\cite{S05,B13} (neglecting the curvature terms) of any unbiased estimator $\widehat{\theta}$ of $\theta=\pi(\overbar{\theta})$ in $\mathcal{M}_{p,k}$ for the proposed error measures~\eqref{eq:icrb:err_hlr} and~\eqref{eq:icrb:err_Gr}.
First of all, exploiting inequality $d_{\mathcal{M}_{p,k}}(\theta,\widehat{\theta})\geq\delta_{\mathcal{M}_{p,k}}^2(\theta,\widehat{\theta})$, we obtain the bound
\begin{equation}
	\mathbb{E}\left[\err^{\mathcal{M}_{p,k}}_{\theta}(\widehat{\theta})\right] \geq \tr\left(\MAT{F}_{\overbar{\theta}}^{\dagger}\right),
\label{eq:icrb:bound_hlr}
\end{equation}
where $\err^{\mathcal{M}_{p,k}}_{\theta}(\widehat{\theta})$ and $\MAT{F}_{\overbar{\theta}}$ are defined in~\eqref{eq:icrb:err_hlr} and~\eqref{eq:icrb:FIM}.

However, as $\MAT{F}_{\overbar{\theta}}$ is adapted to the Riemannian distance $\delta_{\mathcal{M}_{p,k}}$, one cannot expect $\tr(\MAT{F}_{\overbar{\theta}}^{\dagger})$ to well represent the optimal attainable performance when the error is measured with the divergence $d_{\mathcal{M}_{p,k}}$ of proposition~\ref{prop:icrb:div}.
Here, we also conjecture that, since the Fisher information matrix $\MAT{\widetilde{F}}_{\overbar{\theta}}$ defined in~\eqref{eq:icrb:FIM_alt} is constructed from an orthonormal basis on the horizontal space $\widetilde{\mathcal{H}}_{\overbar{\theta}}$ given in~\eqref{eq:icrb:horizontal_alt} which yields the divergence $d_{\mathcal{M}_{p,k}}$, we have the inequality
\begin{equation}
	\mathbb{E}\left[\err^{\mathcal{M}_{p,k}}_{\theta}(\widehat{\theta})\right] \geq \tr\left(\MAT{\widetilde{F}}_{\overbar{\theta}}^{-1}\right).
\label{eq:icrb:bound_hlr_alt}
\end{equation}

Moreover, thanks to the structure of $\MAT{F}_{\overbar{\theta}}$ (see proposition~\ref{prop:icrb:FIM_struct}), it is possible to bound the subspace estimation error~\eqref{eq:icrb:err_Gr}.
Indeed, the block $\MAT{F}_{\MAT{U}_{\perp}}$ is isolated from the rest.
As it is constructed from the elements of the orthonormal basis on $T_{\overbar{\theta}}\overbar{\mathcal{M}}_{p,k}$ of proposition~\ref{prop:icrb:basis} which coincide with the ones of an orthonormal basis on the Grassmann manifold $\mathcal{G}_{p,k}$ associated with the Riemannian distance function $\delta_{\mathcal{G}_{p,k}}$, we have the bound%
\footnote{
	In~\cite{S05}, the only error measure considered is the subspace one~\eqref{eq:icrb:err_Gr}.
	Thus, when dealing with the Gaussian distribution ($\alpha^{++}=1$ in corollary~\ref{corol:icrb:fim}), the bound~\eqref{eq:icrb:bound_Gr} corresponds to the one proposed in~\cite{S05}.
}
\begin{equation}
	\mathbb{E}\left[\err^{\mathcal{G}_{p,k}}_{\theta}(\widehat{\theta})\right] \geq \tr\left(\MAT{F}_{\MAT{U}_{\perp}}^{-1}\right).
\label{eq:icrb:bound_Gr}
\end{equation}
After some manipulations with the basis from proposition \ref{prop:icrb:basis}, it is possible to show that this bound admits the closed-form expression
\begin{equation}
\tr ( \MAT{F}_{\MAT{U}_\bot}^{-1} )
= \frac{(p-k)}{n \alpha^{++}} \sum_{i=1}^k \frac{1 + \sigma_i}{\sigma_i^2},
\end{equation}
where $\{\sigma_i\}_{i=1}^k$ is the set of eigenvalues of $\MAT{\Sigma}$.
As for the Gaussian signal case studied in \cite{S05} (that coincides for $\alpha^{++}=1$ and $\sigma_i = {\rm SNR}, \forall i\in [\![1,k ]\!]$), this leads to an interpretable result in terms of problem dimensions and signal to noise ratio.

\section{Numerical experiments}
\label{sec:num_exp}
\begin{figure*}[t]
\centering
	\setlength\height{5cm} 
\setlength\width{0.5\linewidth}

\begin{tikzpicture}
    \begin{semilogxaxis}[
        width  =\width,
        height =\height,
        at     ={(0,0)},
        %
        xmin        = 10,
        xmax        = 350,
        xtick       = {10,16,100},
        xticklabels = {\empty},
        minor xtick = {20,30,40,50,60,80,90,200,300},
        ylabel = {$\err^{\mathcal{M}_{p,k}}_{\theta}(\widehat{\theta})$~(dB)},
        ymin   = -15,
        ymax   = 10,
        title = {\small$d=3$},
        legend style={legend cell align=left,align=left,draw=none,fill=none,font=\scriptsize,legend columns=2,transpose legend}
        ]
        \addplot[samples=50, smooth,domain=0:6,color=black,dotted,line width=0.8pt,forget plot] coordinates {(16,-15)(16,10)};
        
        \addplot[color=black!50,dashed,line width=0.9pt] table [x=n,y=iCRB_hlr_d3,col sep=comma] {iCRB_hlr_k4.txt};
        \addlegendentry{$\tr(\MAT{F}_{\overbar{\theta}}^{\dagger})$};
        \addplot[color=black,dashed,line width=0.9pt] table [x=n,y=iCRB_hlr_alt_d3,col sep=comma] {iCRB_hlr_k4.txt};
        \addlegendentry{$\tr(\MAT{\widetilde{F}}_{\overbar{\theta}}^{-1})$};

        \addplot[color=black!30,line width=0.9pt,mark=square,mark size=1.6pt] table [x=n,y=pSCM_d3,col sep=comma] {err_hlr_k4.txt};
        \addlegendentry{pSCM};
        \addplot[color=black!50,line width=0.9pt,mark=triangle*, mark options={solid, fill opacity=0, rotate=180},mark size=2.2pt] table [x=n,y=TlrMM_d3,col sep=comma] {err_hlr_k4.txt};
        \addlegendentry{T-MM};
        \addplot[color=black!70,line width=0.9pt,mark=o,mark size=1.6pt] table [x=n,y=TlrRO-st_d3,col sep=comma] {err_hlr_k4.txt};
        \addlegendentry{T-RGD};
        \addplot[color=black!90,line width=0.9pt,mark=triangle*, mark options={solid, fill opacity=0},mark size=2.2pt] table [x=n,y=TlrRO-tr_d3,col sep=comma] {err_hlr_k4.txt};
        \addlegendentry{T-RTR};
    \end{semilogxaxis}

    \begin{semilogxaxis}[
        width  =\width,
        height =\height,
        at     ={(0.84\width,0)},
        %
        xmin        = 10,
        xmax        = 350,
        xtick       = {10,16,100},
        xticklabels = {\empty},
        minor xtick = {20,30,40,50,60,80,90,200,300},
        ymin   = -15,
        ymax   = 10,
        yticklabels  = {\empty},
        title = {\small$d=100$},
        ]
        \addplot[samples=50, smooth,domain=0:6,color=black,dotted,line width=0.8pt,forget plot] coordinates {(16,-15)(16,10)};
        
        \addplot[color=black!50,dashed,line width=0.9pt] table [x=n,y=iCRB_hlr_d100,col sep=comma] {iCRB_hlr_k4.txt};
        \addplot[color=black,dashed,line width=0.9pt] table [x=n,y=iCRB_hlr_alt_d100,col sep=comma] {iCRB_hlr_k4.txt};

        \addplot[color=black!30,line width=0.9pt,mark=square,mark size=1.6pt] table [x=n,y=pSCM_d100,col sep=comma] {err_hlr_k4.txt};
        \addplot[color=black!50,line width=0.9pt,mark=triangle*, mark options={solid, fill opacity=0, rotate=180},mark size=2.2pt] table [x=n,y=TlrMM_d100,col sep=comma] {err_hlr_k4.txt};
        \addplot[color=black!70,line width=0.9pt,mark=o,mark size=1.6pt] table [x=n,y=TlrRO-st_d100,col sep=comma] {err_hlr_k4.txt};
        \addplot[color=black!90,line width=0.9pt,mark=triangle*, mark options={solid, fill opacity=0},mark size=2.2pt] table [x=n,y=TlrRO-tr_d100,col sep=comma] {err_hlr_k4.txt};
    \end{semilogxaxis}

    \begin{semilogxaxis}[
        width  =\width,
        height =\height,
        at     ={(0,-0.71\height)},
        xlabel      = {$n$},
        xmin        = 10,
        xmax        = 350,
        xtick       = {10,16,100},
        xticklabels = {$10^1$,$n=p$,$10^2$},
        minor xtick = {20,30,40,50,60,80,90,200,300},
        ylabel = {$\err^{\mathcal{G}_{p,k}}_{\theta}(\widehat{\theta})$~(dB)},
        ymin   = -22,
        ymax   = 2,
        legend style={legend cell align=left,align=left,draw=none,fill=none,font=\scriptsize,legend columns=2,transpose legend}
        ]
        \addplot[samples=50, smooth,domain=0:6,color=black,dotted,line width=0.8pt,forget plot] coordinates {(16,-22)(16,2)};

        \addplot[color=black,dashed,line width=0.9pt] table [x=n,y=iCRB_Gr_d3,col sep=comma] {iCRB_hlr_k4.txt};
        \addlegendentry{$\tr(\MAT{F}_{\MAT{U}_{\perp}}^{-1})$};
        
        \addlegendimage{empty legend};
        \addlegendentry{};


        \addplot[color=black!30,line width=0.9pt,mark=square,mark size=1.6pt] table [x=n,y=pSCM_d3,col sep=comma] {err_Gr_k4.txt};
        \addlegendentry{pSCM};
        \addplot[color=black!50,line width=0.9pt,mark=triangle*, mark options={solid, fill opacity=0, rotate=180},mark size=2.2pt] table [x=n,y=TlrMM_d3,col sep=comma] {err_Gr_k4.txt};
        \addlegendentry{T-MM};
        \addplot[color=black!70,line width=0.9pt,mark=o,mark size=1.6pt] table [x=n,y=TlrRO-st_d3,col sep=comma] {err_Gr_k4.txt};
        \addlegendentry{T-RGD};
        \addplot[color=black!90,line width=0.9pt,mark=triangle*, mark options={solid, fill opacity=0},mark size=2.2pt] table [x=n,y=TlrRO-tr_d3,col sep=comma] {err_Gr_k4.txt};
        \addlegendentry{T-RTR};
    \end{semilogxaxis}

    \begin{semilogxaxis}[
        width  =\width,
        height =\height,
        at     ={(0.84\width,-0.71\height)},
        xlabel      = {$n$},
        xmin        = 10,
        xmax        = 350,
        xtick       = {10,16,100},
        xticklabels = {$10^1$,$n=p$,$10^2$},
        minor xtick = {20,30,40,50,60,80,90,200,300},
        ymin   = -22,
        ymax   = 2,
        yticklabels  = {\empty},
        ]
        \addplot[samples=50, smooth,domain=0:6,color=black,dotted,line width=0.8pt,forget plot] coordinates {(16,-22)(16,2)};
        
        \addplot[color=black,dashed,line width=0.9pt] table [x=n,y=iCRB_Gr_d100,col sep=comma] {iCRB_hlr_k4.txt};

        \addplot[color=black!30,line width=0.9pt,mark=square,mark size=1.6pt] table [x=n,y=pSCM_d100,col sep=comma] {err_Gr_k4.txt};
        \addplot[color=black!50,line width=0.9pt,mark=triangle*, mark options={solid, fill opacity=0, rotate=180},mark size=2.2pt] table [x=n,y=TlrMM_d100,col sep=comma] {err_Gr_k4.txt};
        \addplot[color=black!70,line width=0.9pt,mark=o,mark size=1.6pt] table [x=n,y=TlrRO-st_d100,col sep=comma] {err_Gr_k4.txt};
        \addplot[color=black!90,line width=0.9pt,mark=triangle*, mark options={solid, fill opacity=0},mark size=2.2pt] table [x=n,y=TlrRO-tr_d100,col sep=comma] {err_Gr_k4.txt};
    \end{semilogxaxis}

\end{tikzpicture}
	\caption{
		Mean of error measures~\eqref{eq:icrb:err_hlr}~(top) and~\eqref{eq:icrb:err_Gr}~(bottom) of methods pSCM, T-MM, T-RGD and T-RTR along with the corresponding intrinsic Cram\'er-Rao bounds~\eqref{eq:icrb:bound_hlr}, \eqref{eq:icrb:bound_hlr_alt} and~\eqref{eq:icrb:bound_Gr} as functions of the number of samples $n$.
		The means are computed over $500$ simulated sets $\{\VEC{x}_i\}_{i=1}^n$ with $d=3$~(left) and $100$~(right), $p=16$ and $k=4$.
	}
\label{fig:sim_k4}
\end{figure*}
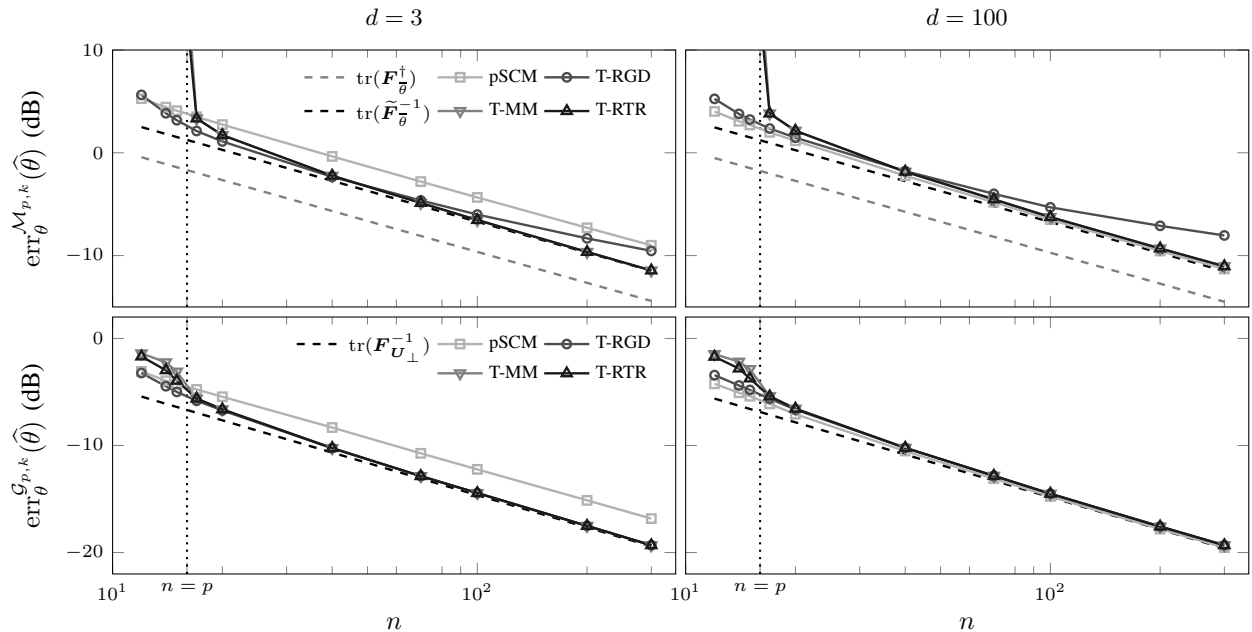

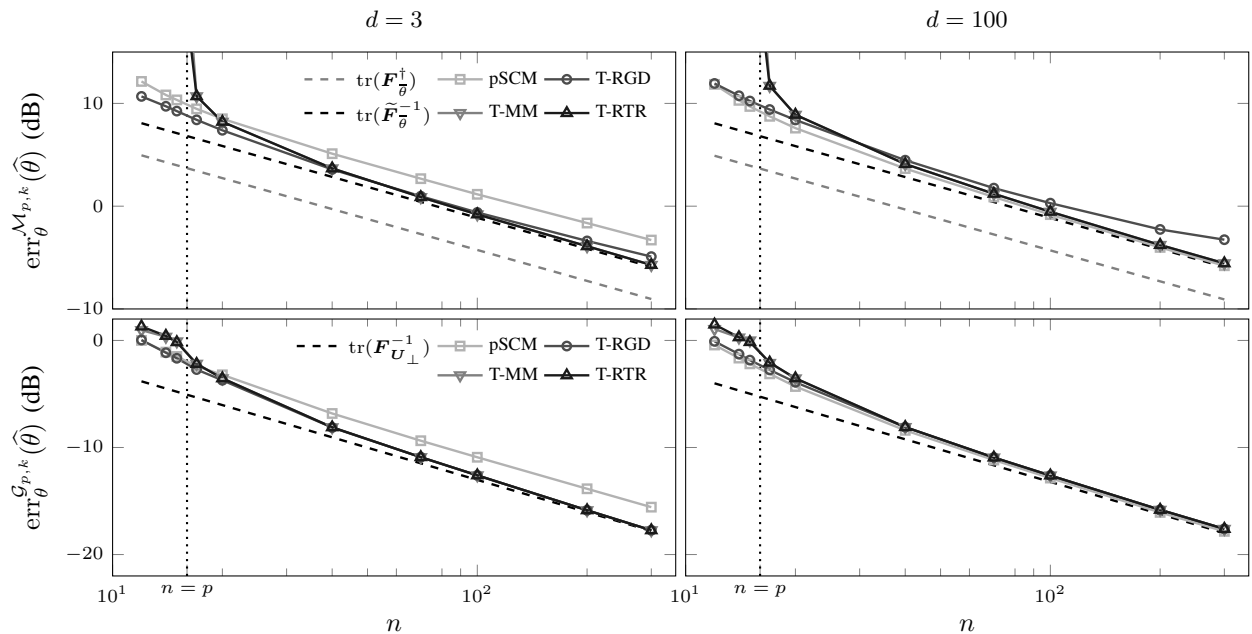
\begin{figure*}[t]
\centering
	\setlength\height{5cm} 
\setlength\width{0.5\linewidth}

\begin{tikzpicture}
    \begin{semilogxaxis}[
        width  =\width,
        height =\height,
        at     ={(0,0)},
        %
        xmin        = 10,
        xmax        = 350,
        xtick       = {10,16,100},
        xticklabels = {\empty},
        minor xtick = {20,30,40,50,60,80,90,200,300},
        ylabel = {$\err^{\mathcal{M}_{p,k}}_{\theta}(\widehat{\theta})$~(dB)},
        ymin   = -10,
        ymax   = 15,
        title = {\small$d=3$},
        legend style={legend cell align=left,align=left,draw=none,fill=none,font=\scriptsize,legend columns=2,transpose legend}
        ]
        \addplot[samples=50, smooth,domain=0:6,color=black,dotted,line width=0.8pt,forget plot] coordinates {(16,-10)(16,15)};
        
        \addplot[color=black!50,dashed,line width=0.9pt] table [x=n,y=iCRB_hlr_d3,col sep=comma] {iCRB_hlr_k8.txt};
        \addlegendentry{$\tr(\MAT{F}_{\overbar{\theta}}^{\dagger})$};
        \addplot[color=black,dashed,line width=0.9pt] table [x=n,y=iCRB_hlr_alt_d3,col sep=comma] {iCRB_hlr_k8.txt};
        \addlegendentry{$\tr(\MAT{\widetilde{F}}_{\overbar{\theta}}^{-1})$};

        \addplot[color=black!30,line width=0.9pt,mark=square,mark size=1.6pt] table [x=n,y=pSCM_d3,col sep=comma] {err_hlr_k8.txt};
        \addlegendentry{pSCM};
        \addplot[color=black!50,line width=0.9pt,mark=triangle*, mark options={solid, fill opacity=0, rotate=180},mark size=2.2pt] table [x=n,y=TlrMM_d3,col sep=comma] {err_hlr_k8.txt};
        \addlegendentry{T-MM};
        \addplot[color=black!70,line width=0.9pt,mark=o,mark size=1.6pt] table [x=n,y=TlrRO-st_d3,col sep=comma] {err_hlr_k8.txt};
        \addlegendentry{T-RGD};
        \addplot[color=black!90,line width=0.9pt,mark=triangle*, mark options={solid, fill opacity=0},mark size=2.2pt] table [x=n,y=TlrRO-tr_d3,col sep=comma] {err_hlr_k8.txt};
        \addlegendentry{T-RTR};
    \end{semilogxaxis}

    \begin{semilogxaxis}[
        width  =\width,
        height =\height,
        at     ={(0.84\width,0)},
        %
        xmin        = 10,
        xmax        = 350,
        xtick       = {10,16,100},
        xticklabels = {\empty},
        minor xtick = {20,30,40,50,60,80,90,200,300},
        ymin   = -10,
        ymax   = 15,
        yticklabels  = {\empty},
        title = {\small$d=100$},
        ]
        \addplot[samples=50, smooth,domain=0:6,color=black,dotted,line width=0.8pt,forget plot] coordinates {(16,-10)(16,15)};
        
        \addplot[color=black!50,dashed,line width=0.9pt] table [x=n,y=iCRB_hlr_d100,col sep=comma] {iCRB_hlr_k8.txt};
        \addplot[color=black,dashed,line width=0.9pt] table [x=n,y=iCRB_hlr_alt_d100,col sep=comma] {iCRB_hlr_k8.txt};

        \addplot[color=black!30,line width=0.9pt,mark=square,mark size=1.6pt] table [x=n,y=pSCM_d100,col sep=comma] {err_hlr_k8.txt};
        \addplot[color=black!50,line width=0.9pt,mark=triangle*, mark options={solid, fill opacity=0, rotate=180},mark size=2.2pt] table [x=n,y=TlrMM_d100,col sep=comma] {err_hlr_k8.txt};
        \addplot[color=black!70,line width=0.9pt,mark=o,mark size=1.6pt] table [x=n,y=TlrRO-st_d100,col sep=comma] {err_hlr_k8.txt};
        \addplot[color=black!90,line width=0.9pt,mark=triangle*, mark options={solid, fill opacity=0},mark size=2.2pt] table [x=n,y=TlrRO-tr_d100,col sep=comma] {err_hlr_k8.txt};
    \end{semilogxaxis}

    \begin{semilogxaxis}[
        width  =\width,
        height =\height,
        at     ={(0,-0.71\height)},
        xlabel      = {$n$},
        xmin        = 10,
        xmax        = 350,
        xtick       = {10,16,100},
        xticklabels = {$10^1$,$n=p$,$10^2$},
        minor xtick = {20,30,40,50,60,80,90,200,300},
        ylabel = {$\err^{\mathcal{G}_{p,k}}_{\theta}(\widehat{\theta})$~(dB)},
        ymin   = -22,
        ymax   = 2,
        legend style={legend cell align=left,align=left,draw=none,fill=none,font=\scriptsize,legend columns=2,transpose legend}
        ]
        \addplot[samples=50, smooth,domain=0:6,color=black,dotted,line width=0.8pt,forget plot] coordinates {(16,-22)(16,2)};

        \addplot[color=black,dashed,line width=0.9pt] table [x=n,y=iCRB_Gr_d3,col sep=comma] {iCRB_hlr_k8.txt};
        \addlegendentry{$\tr(\MAT{F}_{\MAT{U}_{\perp}}^{-1})$};
        
        \addlegendimage{empty legend};
        \addlegendentry{};


        \addplot[color=black!30,line width=0.9pt,mark=square,mark size=1.6pt] table [x=n,y=pSCM_d3,col sep=comma] {err_Gr_k8.txt};
        \addlegendentry{pSCM};
        \addplot[color=black!50,line width=0.9pt,mark=triangle*, mark options={solid, fill opacity=0, rotate=180},mark size=2.2pt] table [x=n,y=TlrMM_d3,col sep=comma] {err_Gr_k8.txt};
        \addlegendentry{T-MM};
        \addplot[color=black!70,line width=0.9pt,mark=o,mark size=1.6pt] table [x=n,y=TlrRO-st_d3,col sep=comma] {err_Gr_k8.txt};
        \addlegendentry{T-RGD};
        \addplot[color=black!90,line width=0.9pt,mark=triangle*, mark options={solid, fill opacity=0},mark size=2.2pt] table [x=n,y=TlrRO-tr_d3,col sep=comma] {err_Gr_k8.txt};
        \addlegendentry{T-RTR};
    \end{semilogxaxis}

    \begin{semilogxaxis}[
        width  =\width,
        height =\height,
        at     ={(0.84\width,-0.71\height)},
        xlabel      = {$n$},
        xmin        = 10,
        xmax        = 350,
        xtick       = {10,16,100},
        xticklabels = {$10^1$,$n=p$,$10^2$},
        minor xtick = {20,30,40,50,60,80,90,200,300},
        ymin   = -22,
        ymax   = 2,
        yticklabels  = {\empty},
        ]
        \addplot[samples=50, smooth,domain=0:6,color=black,dotted,line width=0.8pt,forget plot] coordinates {(16,-22)(16,2)};
        
        \addplot[color=black,dashed,line width=0.9pt] table [x=n,y=iCRB_Gr_d100,col sep=comma] {iCRB_hlr_k8.txt};

        \addplot[color=black!30,line width=0.9pt,mark=square,mark size=1.6pt] table [x=n,y=pSCM_d100,col sep=comma] {err_Gr_k8.txt};
        \addplot[color=black!50,line width=0.9pt,mark=triangle*, mark options={solid, fill opacity=0, rotate=180},mark size=2.2pt] table [x=n,y=TlrMM_d100,col sep=comma] {err_Gr_k8.txt};
        \addplot[color=black!70,line width=0.9pt,mark=o,mark size=1.6pt] table [x=n,y=TlrRO-st_d100,col sep=comma] {err_Gr_k8.txt};
        \addplot[color=black!90,line width=0.9pt,mark=triangle*, mark options={solid, fill opacity=0},mark size=2.2pt] table [x=n,y=TlrRO-tr_d100,col sep=comma] {err_Gr_k8.txt};
    \end{semilogxaxis}

\end{tikzpicture}
	\caption{
		Mean of error measures~\eqref{eq:icrb:err_hlr}~(top) and~\eqref{eq:icrb:err_Gr}~(bottom) of methods pSCM, T-MM, T-RGD and T-RTR along with the corresponding intrinsic Cram\'er-Rao bounds~\eqref{eq:icrb:bound_hlr}, \eqref{eq:icrb:bound_hlr_alt} and~\eqref{eq:icrb:bound_Gr} as functions of the number of samples $n$.
		The means are computed over $500$ simulated sets $\{\VEC{x}_i\}_{i=1}^n$ with $d=3$~(left) and $100$~(right), $p=16$ and $k=8$.
	}
\label{fig:sim_k8}
\end{figure*}

This section illustrates our Riemannian optimization framework and performance analysis for robust covariance estimation.
In order to do so, we perform covariance estimation of simulated data drawn from the multivariate Student $t$-distribution with $d=3$ (highly non-Gaussian) and $d=100$ (almost Gaussian) degrees of freedom; see~\cite{OTKP12a} for details.

To generate a covariance matrix admitting the structure~\eqref{eq:model:spiked_model}, we compute
\begin{equation*}
	\MAT{R} = \MAT{I}_p + \sigma\MAT{U}\MAT{\Sigma}\MAT{U}^H,
\end{equation*}
where 
\begin{itemize}
	\item $\MAT{U}$ is a random matrix in $\textup{St}_{p,k}$,
	\item $\MAT{\Sigma}$ is a diagonal matrix whose minimal and maximal elements are $\nicefrac1{\sqrt{c}}$ and $\sqrt{c}$ ($c=20$ is the condition number with respect to inversion of $\MAT{\Sigma}$); its other elements are randomly drawn from the uniform distribution between $\nicefrac1{\sqrt{c}}$ and $\sqrt{c}$; its trace is then normalized as $\tr(\MAT{\Sigma})=\tr(\MAT{I}_k)=k$,
	\item $\sigma=50$ is a free parameter corresponding to the spike to noise ratio.
\end{itemize}
In our experiment, we choose $p=16$ and $k\in\{4,8\}$.
Sets $\{\VEC{x}_i\}_{i=1}^n$ are drawn from the multivariate Student $t$-distribution with covariance $\MAT{R}$ and $d\in\{3,100\}$, where $n\in\{12,14,15,17,20,40,70,100,200,300\}$.
For each value of $n$, $500$ sets $\{\VEC{x}_i\}_{i=1}^n$ are simulated and the aim is to estimate the structured covariance matrix $\MAT{R}$ in each~case.

The considered estimators in this experiment are:
\begin{enumerate}[label=(\alph*)]
	\item Projected sample covariance matrix $\MAT{I}_p+\varphi(\widehat{\theta}_{\textup{pSCM}})$ obtained by projecting $n^{-1}\sum_i \VEC{x}_i\VEC{x}_i^H$ on $\MAT{I}_p+\mathcal{H}^{+}_{p,k}$ with~\cite[equation (53)]{SBP16}.
	\item Structured Tyler's $M$-estimator $\MAT{I}_p+\varphi(\widehat{\theta}_{\textup{T-MM}})$ solved with~\cite[algorithm 5]{SBP16}.
	\item Structured Tyler's $M$-estimator $\MAT{I}_p+\varphi(\widehat{\theta}_{\textup{T-RGD}})$ solved with a Riemannian gradient descent algorithm on $\mathcal{M}_{p,k}$; see~\cite[chapter 4]{AMS08}.
	\item Structured Tyler's $M$-estimator $\MAT{I}_p+\varphi(\widehat{\theta}_{\textup{T-RTR}})$ solved with a Riemannian trust region algorithm (second order optimization method) on $\mathcal{M}_{p,k}$; see~\cite[chapter 7]{AMS08}.
\end{enumerate}
The three iterative methods are initialized with the principal subspace of the projected sample covariance matrix estimator, \emph{i.e.}, $(\MAT{\widehat{U}}_{\textup{pSCM}},\MAT{I}_k)$.
Riemannian optimization on $\mathcal{M}_{p,k}$ is performed with manopt toolbox~\cite{BMAS14} and we choose $\alpha=\frac{p+d}{p+d+1}$ and $\beta=\alpha-1$ in the Riemannian metric~\eqref{eq:rg_hlr:metric}.

In figures~\ref{fig:sim_k4} and~\ref{fig:sim_k8}, we observe that, in all considered cases, \emph{i.e.} $d\in\{3,100\}$ and $k\in\{4,8\}$, the lower bound~\eqref{eq:icrb:bound_hlr} is not reached by any of the methods for error measure~\eqref{eq:icrb:err_hlr}.
This is expected as this bound is suited to the Riemannian distance on $\mathcal{M}_{p,k}$ and not to the divergence of proposition~\ref{prop:icrb:div}.
However, for error measure~\eqref{eq:icrb:err_hlr}, the bound~\eqref{eq:icrb:bound_hlr_alt}, which arises from the Fisher information matrix well suited to our divergence, is reached by several methods as the number of samples $n$ grows.
Concerning the subspace error~\eqref{eq:icrb:err_Gr}, the lower bound~\eqref{eq:icrb:bound_Gr} is reached in all considered cases by several methods as $n$ grows.
Further notice that, for $k=4$, a smaller amount of samples $n$ is needed for the bounds~\eqref{eq:icrb:bound_hlr_alt} and~\eqref{eq:icrb:bound_Gr} to be attained than for $k=8$.

Unlike the other considered estimators, the performance of pSCM depends on the degree of freedom $d$ of the Student $t$-distribution.
As expected, when data are close to Gaussianity ($d=100$), pSCM provides good results and attains both bounds~\eqref{eq:icrb:bound_hlr_alt} and~\eqref{eq:icrb:bound_Gr}.
However, when they are far from being Gaussian ($d=3$), pSCM fails to give optimal results.
We also observe that T-MM and T-RTR have very similar performance.
They both fail when $n$ is small, especially when it gets close to $p$ (or smaller).
However, they perform well when $n$ is sufficient and reach both bounds~\eqref{eq:icrb:bound_hlr_alt} and~\eqref{eq:icrb:bound_Gr}.
Concerning T-RGD, we notice that it yields good results as compared to other estimators when $n$ is small.
As $n$ grows, even though T-RGD still provide satisfying subspaces (bound~\eqref{eq:icrb:bound_Gr} is reached by error measure~\eqref{eq:icrb:err_Gr}), its performance with respect to error measure~\eqref{eq:icrb:err_hlr} deteriorates as compared to other estimators.
In conclusion, our optimization framework on $\mathcal{M}_{p,k}$ provides satisfying results on these simulated data for all considered cases.
Depending on the number of samples at hand, different optimization algorithms are preferable: the first order method (T-RGD) is more advantageous when a small amount of samples is available whereas the second order method (T-RTR) performs better as the number of samples grows.

\section{Conclusions and perspectives}
\label{sec:conclusion}
This article proposes an original Riemmanian geometry to study low-rank structured elliptical models.
The tools developed within this framework (representations of tangent spaces, geodesics, Riemannian gradient and Hessian, retraction, divergence function) allow to derive both estimation algorithms and intrinsic Cram\'er-Rao lower bounds adapted to these models with a unified view.
Some potential extensions of this work include: generalization to $M$-estimators and estimation of the parameters of the Fisher information metric, integration of curvature terms and intrinsic bias in the intrinsic Cram\'er-Rao lower bounds.

\ifCLASSOPTIONcaptionsoff
  \newpage
\fi



%
\bibliographystyle{unsrt}
\bibliography{biblio}

\end{document}